\documentclass[11pt,a4paper,oneside,english,parskip=full,abstracton]{scrartcl}

\usepackage[top=2cm, bottom=2cm, left=2cm, right=1.5cm]{geometry}

\usepackage[utf8]{inputenc}
\usepackage[T1]{fontenc}
\usepackage[english]{babel}
\usepackage{graphicx}
\usepackage{csquotes}
\usepackage{amsmath}
\usepackage{amssymb}
\usepackage{amsthm}
\usepackage{algorithm}
\usepackage{algpseudocode}
\usepackage{authblk}
\usepackage{subcaption}
\usepackage{caption}
\usepackage[capitalise,noabbrev]{cleveref}
\usepackage{xstring}
\usepackage[bottom]{footmisc}

\usepackage[
backend=biber,
citestyle=numeric-comp,
firstinits=true
]{biblatex}

\usepackage{tikz}
\usepackage{pgfplots}
\usepackage{pgfplotstable,booktabs}
\usetikzlibrary{calc,math}

\theoremstyle{plain}
\newtheorem{theorem}{Theorem}[section] 
\theoremstyle{definition}

\newtheorem{remark}[theorem]{Remark}
\newtheorem{problem}[theorem]{Problem}

\newcommand{\ignore}[1]{}

\begin{document}
	
	\selectlanguage{english}
	\pagenumbering{arabic}
	
	\title{Matrix-free multigrid solvers\\ for phase-field fracture problems}
	
	\author[1]{D. Jodlbauer}
	\author[1]{U. Langer}
	\author[2,3]{T. Wick}
	
	\affil[1]{Johann Radon Institute for Computational and Applied Mathematics, Austrian Academy of Sciences, Altenbergerstr. 69, A-4040 Linz, Austria}
	\affil[2]{Institut f\"ur Angewandte Mathematik, Leibniz Universit\"at Hannover, Welfengarten 1, 30167 Hannover, Germany}
	\affil[3]{Cluster of Excellence PhoenixD (Photonics, Optics, and
		Engineering - Innovation Across Disciplines), Leibniz Universit\"at Hannover, Germany}
	
	\date{}
	
	\maketitle
	
	
	\begin{abstract}	
		In this work, we present a framework for the matrix-free solution to a 
		monolithic quasi-static phase-field fracture model with geometric multigrid methods.
		Using a standard matrix based approach 
		within the Finite Element Method
		requires lots of memory, 
		which eventually becomes a serious bottleneck.
		A matrix-free approach overcomes this problems and greatly reduces the 
		amount of required memory, allowing to solve larger problems on available
		hardware. One key challenge is concerned with the crack irreversibility 
		for which a primal-dual active set method is employed. Here, the active 
		set values of fine meshes must be available on coarser levels of the multigrid
		algorithm. The developed multigrid method 
		provides
		a preconditioner for 
		a generalized minimal residual (GMRES) solver. This method is used 
		for solving the linear equations inside Newton's method for treating the 
		overall nonlinear-monolithic discrete displacement/phase-field formulation.
		Several numerical examples demonstrate 
		the performance and robustness of our solution technology.
		Mesh refinement studies, variations in the phase-field 
		regularization parameter, 
		iterations numbers of the linear and nonlinear solvers, and some parallel 
		performances are conducted to substantiate the efficiency of the proposed
		solver for single fractures, multiple pressurized fractures, and 
		a L-shaped panel test in three dimensions. \\[1em]
		
	\end{abstract}


	\section{Introduction}
	
	Predicting fracture growth can be of big interest in a variety of fields, such as manufacturing processes, engineering sciences, and medical applications.
	However, fracture propagation poses several challenges for numerical methods.
	A main problem is the fact, that cracks are usually lower-dimensional phenomena.
	For example, fractures in a sheet of glass are, most often, just thin lines.
	This makes it very difficult for standard numerical methods, 
	like the Finite Element Method (FEM)
	to represent such kind of discontinuities.
	Another challenge is the irreversibility property, i.e. a fracture
	should not be able to heal itself, which introduces time-dependent
	constraints on the solution. More difficulties are related to 
	fracture nucleation, branching, and fracture networks with
	curve-linear, complex, crack patterns.

	Methods like generalized/extended FEM (GFEM/XFEM) try to overcome these fallacies of standard FEM by enriching the solution space with discontinuous functions.
	This allows the representation of sharper objects like cracks, see
	e.g. \cite{BeBl99,MeDu07}. A short overview of various methods is compiled in \cite{WiSiWh16}.
	
	We now focus on a variational phase-field  based approach,
	which was first introduced for fractures in
	\cite{BoFrMa00} based on the fracture model presented in \cite{FrMa98}.
	Therein, Griffith's original model for brittle fractures \cite{GrTa21} is rewritten by means of an energy minimization problem.
	This numerical approximation introduces an additional scalar-valued indicator variable representing the fracture.
	It is then computed as the minimizer of the total energy functional, together with the unknown solid displacement.
	Specifically therein, the continuous model is approximated using elliptic Ambrosio-Tortorelli functionals, see \cite{AmTo90}.
	Such an approximation yields a smooth representation of the originally sharp crack.
	Thus, the solution does no longer jump from broken to unbroken parts, but provides a smoothed transition between them.
	The approximation quality is given by the phase-field regularization parameter $\varepsilon$.
	Roughly speaking, it determines the size of the diffusive zone around the fracture, see \cref{pic:notation} for an illustration.
	
	This smoothing effect could be seen as a disadvantage of phase-field
	methods, as it is no longer possible to accurately localize the
	fracture. Indeed, fracture width computations \cite{LeWhWi17} and 
	posing interface conditions on the fracture boundary are challenging 
	\cite{MiWhWi15} with no final answer to date.
	On the other hand, it comes with some benefits: the problem is reduced
	to find the minimizers of a nonlinear energy functional, which is a
	common procedure in engineering and applied mathematics.
	Thus, it is possible to use well-known Galerkin techniques, 
	making it easier to find suitable software to assist with the
	implementation. Furthermore, the fracture is entirely represented using an additional variable.
	In particular, the mesh is fixed and does not need to be updated when 
	the fracture propagates, as it is the case for methods like XFEM.

	The phase-field fracture (PFF) model can be enhanced to include, for
	instance, 
	pressure-driven cracks; see, e.g,
	\cite{MiMaTe15,MiMa16,LeWhWi16,MiWhWi15a,MiWhWi19,HeMa17,WiLa16,YoBo16,ChBoYo19,SiVe18}. 
	Furthermore, thermodynamically consistent schemes have been developed in \cite{MiHoWe10,MiWeHo10}.
	Therein, the elastic energy is split into tensile and compressive parts.
	Other splitting strategies exists, e.g. \cite{AmMaMa09} based on
	a deviatoric split, or a hybrid strategy as presented in
	\cite{AmGeDe15}. The latter also contains a nice comparison of the different splitting techniques.
	
	The most challenging part is related to the nonlinear equations and 
	the irreversibility condition, which appears in terms of a variational inequality.
	For treating the irreversibility constraint, several approaches
	(simple penalization, augmented Lagrangian techniques, or a
	strain-history field) have
	been proposed in the literature. In this work, we adopt 
	a primal-dual active set method \cite{HiItKu03,ItKu00} that 
	can be identified as a semi-smooth Newton method \cite{HiItKu03} and 
	was first applied to phase-field fracture in \cite{HeWhWi15}.
	With regard to the nonlinearities in the governing equations
	for displacements and phase-field two main approaches can be
	distinguished: staggered approaches \cite{Bo07,BuOrSu10,MiHoWe10}
	and monolithic methods \cite{GeDe16,Wi17,Wi17b}.
	
	Typically, the nonlinear solvers reduce the problem to the repeated solution of linear systems of equations.
	Krylov-subspace based solvers like Conjugate Gradients (CG) or Generalized Minimal Residual (GMRES) \cite{Sa03} have been proven to be very effective methods for solving such large sparse linear systems, if provided with a reasonable preconditioner.
	Different possibilities to precondition the PFF problem are shown in \cite{FaMa17}.
	
	A powerful class of solvers/preconditioners that is known to work for many different types of problems are multigrid methods. 
	These methods smooth the iteration error on the fine grid and correct the iterate by 
	a course-grid correction in a two-grid setting, whereas a recurrent applications on a sequence of grids
	lead to a multigrid procedure; see  \cite{Ha85,Br93}.
	Multigrid methods have been successfully applied to a wide range of related problems, see e.g. 
	\cite{Ko94,Ko02,HaMi83,Ho87,HoMi89,Be93,BlBrSu04,ReVoHe13,HeLaNa13,GrSa19,KoGr09}. 
	Algebraic multigrid methods provided by Trilinos \cite{HeBaHo03,GeSiHu06} have been applied to phase-field fracture in \cite{HeWi18}, yielding a scalable parallel solution scheme; see also 
	\cite{St01,HaLa02,BrScSc07} for algebraic 
	multigrid method and their parallelization.	
	
	When attempting to solve larger and larger problems, memory soon becomes the limiting factor.
	To overcome this, matrix-free methods have been developed.
	Instead of assembling and storing the discrete Jacobian, as it is usually done, the necessary information is computed on-the-fly.
	This drastically reduces the amount of memory required, but may increase the computational cost.
	Hence, a clever implementation is needed to keep up with standard matrix-based methods.
	In this work, we rely on the FEM library deal.II \cite{AlArBa18},
	which readily includes such an implementation 
	\cite{KrKo12,KrKo17}.
	
	As indicated in our previous overviews, the main aims of this paper are two-fold:
	\begin{itemize}
		\item Development of a geometric multigrid preconditioner for the GMRES
		iterations inside Newton's method;
		\item A matrix-free implementation of the geometric multigrid method.
	\end{itemize}
	Here, one focus is on quantities that are given naturally on the
	finest grid and which may also be required on the coarser levels of the multigrid algorithm, e.g. the current linearization point or the Active-Set.
	We will present strategies to handle those problems with relative
	ease. Indeed, the FEM library deal.II has gone to great lengths in
	providing a rich and flexible interface for the implementation of such sophisticated methods.	
	
	The outline of this paper is as follows.
	In \cref{sec:problem}, we will describe the governing equations of phase-field fractures and clarify some of the notation used.
	Furthermore, we briefly describe some details on the general discretization scheme used.
	In \cref{sec:solution}, we discuss our approach for the solution of the phase-field problem.
	This includes the non-linear solver (Active-Set) and the linear solver (geometric Multigrid).
	The main part of this work is presented in \cref{sec:matrixfree}, which is devoted to the matrix-free realization of this solution approach.
	Finally, we present our findings and numerical results for various test cases in \cref{sec:results}.
	
	
	\section{Problem Description}
	\label{sec:problem}
	
	In this section, we introduce the basic notation and the underlying equations.
	In the following, let $D \subset \mathbb{R}^d, d=2,3$ the total domain wherein $\mathcal{C}\subset \mathbb{R}^{d-1}$ denotes the fracture and $\Omega \subset \mathbb{R}^d$ {is} the intact domain.
	We assume (possibly time-dependent non-homogeneous) Dirichlet conditions on the outer boundary $\partial D$.

	\begin{figure}[ht]
		\begin{center}
			\scalebox{0.6}{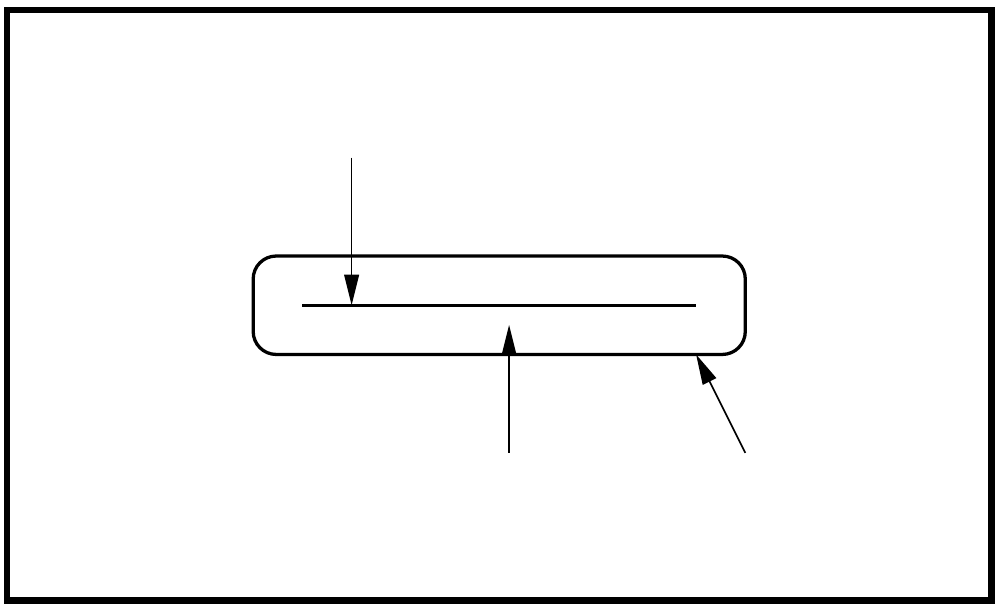}\hspace*{1cm}
			\includegraphics[width=7cm,height=3.75cm]{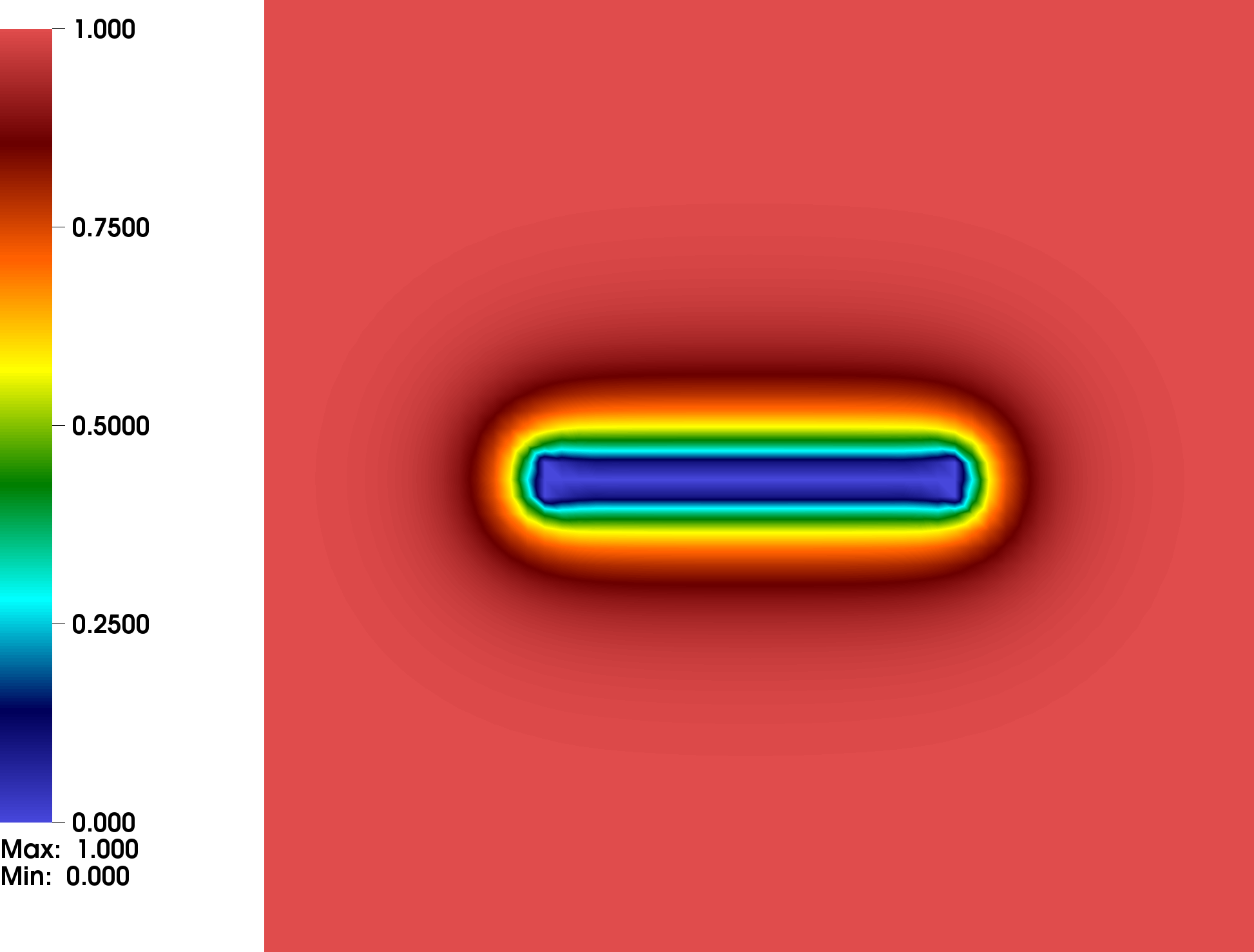}
		\end{center}
		
		\caption{Illustration of the Notation: sharp fracture $\mathcal{C}$ and phase-field approximation $\Omega_F$.
		}
		
		\label{pic:notation}
	\end{figure}
	
	In the phase-field approach, the one-dimensional fracture $\mathcal{C}$ is approximated by $\Omega_F \subset \mathbb{R}^d$ with the help of an elliptic (Ambrosio-Tortorelli) functional \cite{AmTo90,AmTo92}.
	Using the phase-field variable $\varphi$ (introduced in \cref{sec:problem}), $\Omega_F$ can be defined as
	\begin{equation}
	\label{eq:omega}
	\Omega_F := \{ x\in \mathbb{R}^d \vert \; \varphi(x) < 1 \}.
	\end{equation}
	The lower-dimensional fracture $\mathcal{C}$ and its phase-field approximation $\Omega_F$ are visualized in \cref{pic:notation}.
	The width of the phase-field fracture is determined by $\varepsilon$.
	For fracture formulations posed in a variational setting, this has been first proposed in \cite{BoFrMa00}.
	The inner fracture boundary is denoted by the smeared $\varepsilon$-dependent boundary $\partial\Omega_F$.
	We note that the precise location of $\partial\Omega_F$ is of no
	importance in this work 
	(in contrast to e.g., \cite{MiWhWi15}).
	The reader is referred to \cref{pic:notation} for an illustration of the notation.
	Finally, we denote the $L^2$ scalar product with $(a,b) := (a,b)_{D}:=\int_{D} a \cdot b \ dx$.
	For tensor-valued functions $A,B$, we use $(A,B) := (A,B)_{D}:=\int_{D} A : B \ dx$.
	
	\subsection{Quasi-static phase-field for brittle fracture}
	
	We briefly recapitulate the ingredients for a phase-field model for
	solid mechanics and pressurized fractures in brittle materials.
	Such a model is based on the variational/phase-field fracture approach of \cite{FrMa98,BoFrMa00}.
	Thermodynamically-consistent phase-field techniques using a stress-split into tension and compression have been proposed in \cite{AmMaMa09} and \cite{MiWeHo10}.
	
	The previous formulations start with an energy functional $E(u, \varphi)$ which is minimized with respect to the unknown solution variables $u : D \to \mathbb{R}^d$ (displacements) and a smoothed scalar-valued phase-field function $\varphi: D \to [0,1]$.
	The latter one varies in the zone of size $\varepsilon$ from $0$ (fracture) to $1$ (intact material).
	The first-order necessary condition are the Euler-Lagrange equations, which are obtained by differentiation with respect to the two unknowns $u$ and $\varphi$.
	Adding a pressure $p : D \to \mathbb{R}$ to the Euler-Lagrange
	equations that acts on the fracture boundary has been formulated and
	mathematically analyzed in \cite{MiWhWi15a,MiWhWi19} for staggered and
	monolithic approaches, respectively.
	In all the previous fracture models, the physics of the underlying problem ask to enforce a crack irreversibility condition (the crack can never heal) that is an inequality condition in time:
	\begin{align}
	\partial_t \varphi \leq 0.
	\end{align}
	Consequently, modeling of fracture evolution problems leads to a variational inequality system, that is always, due to this constraint, quasi-stationary or time-dependent.
	The resulting variational formulation is stated in an incremental (i.e., time-discretized) formulation in which the continuous irreversibility constraint is approximated by \[ \varphi \leq \varphi^{old}. \]
	Here, $\varphi^{old}$ will later denote the previous time step solution and $\varphi$ the current solution.
	Let $V:=H^1_0(D)$ and 
	\begin{equation}
	\label{W_in}
	W_{in}:=\{w\in H^1(D) |\, w\leq \varphi^{old} \leq 1 \text{ a.e. on
	} D\} 
	\end{equation}
	be the function spaces we work with here; and for later purposes we also need $W:=H^1(D)$.
	The Euler-Lagrange system for pressurized phase-field fractures reads \cite{MiWhWi19}:
	\begin{problem}
		\label{problem}
		Let $p\in L^{\infty}(D)$ be given.
		Find vector-valued displacements and a scalar-valued phase-field variable $\{ u, \varphi \}:= \{ u^n, \varphi^n \} \in \{u_D + V\} \times W_{in}$ such that:
		\begin{align*}
		& \left(  g_\kappa (\varphi) \sigma^+(u), e(w) \right) + \left(\sigma^-(u), e(w) \right) + \left(\varphi^2 p, \mbox{div } w \right) = 0 \quad & \forall w \in V, \\
		& & \\
		& \frac{1}{2} \left(  \partial_\varphi g_\kappa(\varphi) E_s^+(e(u)), \psi - \varphi \right) + 2 (\varphi \ p \ \mbox{div } u, \psi - \varphi) & \\
		& + G_c \left( -\frac{1}{\varepsilon} ( 1-\varphi, \psi -\varphi ) + \varepsilon \left( \nabla \varphi, \nabla (\psi - \varphi ) \right) \right) \geq 0 \quad & \forall \psi \in W \cap L^\infty(D),
		\end{align*}
		at each incremental step $n=1,2,3,\ldots$.
	\end{problem}
	The elastic energy is denoted by $$E_s(e) := \frac{1}{2} \lambda \mbox{tr}(e)^2 + \mu (e, e)$$ with resulting stress tensor $$\sigma := \frac{\partial E_s(e)}{\partial e} = 2\mu e + \lambda \mbox{tr}(e) I.$$
	The parameters $\mu$ and $\lambda$ denote the usual Lam\'e coefficients, $e(u) = \frac{1}{2} (\nabla u + \nabla u^T)$ is the linearized strain tensor and $I$ denotes the identity matrix.
	The degradation function $g_\kappa$ is defined as $g_\kappa(\varphi) := (1-\kappa) \varphi^2 + \kappa$, with a regularization parameter $\kappa \ll 1$.
	Different choices for $g_\kappa$ are possible (see \cite{KuScMl15}), for example $g_\kappa := \varphi^2 + \kappa$ as used in \cite{BoFrMa00}.
	Physically, $\kappa$ represents the residual stiffness of the material.
	Consequently, since \[ g_\kappa(\varphi) \rightarrow \kappa \quad \text{for } \varphi \to 0, \] the material stiffness decreases while approaching the fracture zone.
	
	The pressure terms $( \varphi^2 p, \mbox{div } w)$ and $2 ( \varphi \ p\; \mbox{div } u, \psi - \varphi )$ have been derived in \cite{MiWhWi15a,MiWhWi19} and are based on an interface law that has been further manipulated using Gauss' divergence theorem.
	
	A superscript $+$ denotes the splitting of energy/stress into two parts, i.e. $E_s = E_s^+ + E_s^-$ and consequently, similar for $\sigma$.
	In this work, we mainly considered the case without splitting ($E_s^+ := E_s$) and the Miehe-type splitting $$E_s^+ := \frac{1}{2} \lambda \langle \mbox{tr}(e) \rangle_+^2 + \mu \mbox{tr}(e_+^2); $$ see \cite{MiWeHo10,MiHoWe10}.
	The positive part of a scalar variable is defined as $\langle x \rangle_+ := \mbox{max}(0, x)$.
	For the strain tensor, the positive part is defined by means of the eigenvalues.
	More precisely, $$e_+ := \sum_{i=1}^{d} \langle \lambda_i \rangle_+ n_i \otimes n_i,$$ with $\lambda_i$ and $n_i$ denoting the eigenvalues and eigenvectors of the strain tensor $e$.
	Consequently, the stress tensor splitting is defined as $\sigma^+(u) = \partial_e E_s^+$ (independent of the actual splitting used).
	A discussion of various splitting laws can be found in \cite{BoVeSc12}[Section 2.2] and \cite{AmGeDe15}.
	
	\subsection{Challenges}
	\label{sec:challenges}
	
	In the following, we will briefly describe the main challenges regarding the phase-field fracture propagation problem.
	
	\textbf{Relationship of $\varepsilon - h$}

	From the theory of $\Gamma$-convergence (e.g., \cite{AmTo92}) 
	we want to have $h \in o(\varepsilon)$.
	Intuitivly speaking, we require a mesh-size $h$ small enough to resolve $\varepsilon$.
	Hence, we do not only have to investigate convergence properties with
	respect to $h$, but also to $\varepsilon$.
	Due to the interplay of $h$ and $\epsilon$, this is a serious
	theoretical challenge.
	Some first numerical experiments were undertaken in \cite{Wi16a}.
	
	\textbf{Nonlinearities}
	
	Problem \Cref{problem} is nonlinear due to the monolithic formulation, the stress splitting, and the inequality constraint.
	The most critical term is the quasi-linearity $\left(g_\kappa(\varphi)
	\sigma^+(u), e(w) \right)$, which causes most of the challenges in
	designing reliable and efficient solution algorithms because
	the corresponding Hessian matrix is indefinite. 
	In fact, integrating this term with respect to $u$ yields the corresponding term on the energy level:
	\[ \left( g_\kappa(\varphi) \sigma^+(u), e(u) \right), \]
	which has been well-characterized to be non-convex in both variables
	$u$ and $\varphi$ simultaneously in the very early work
	\cite{BoFrMa00,Bo07} (without the stress-splitting though).
	A very simple 
	prototype example demonstrating the difficulties was studied in \cite{Wi17b}.
	
	\textbf{Linear solver}
	
	Another major issue is the solution of the linearized systems of equations.
	Depending on the linearization strategy used, we will have to solve many linear equations.
	Thus, in order to simulate large problems, a fast and robust linear
	solver is an absolute necessity.
	Here, a crucial issue is the robustness of the linear solver with
	respect to discretization, model, and material parameters.
	
	\textbf{Matrix-Free framework}
	
	In this work, we will consider a matrix-free implementation for the geometric multigrid preconditioner.
	Implementing such a solver is quite challenging.
	For example, quantities given naturally on the finest grid may also be required on the coarser levels of the multigrid algorithm, e.g. the current linearization point or the Active-Set.
	In this paper we will present strategies to handle those problems with relative ease.
	Furthermore, we would like to mention that the FEM library deal.II has gone to great lengths in providing a rich and flexible interface for the implementation of such sophisticated methods.
	
	\textbf{Adaptivity and parallelization}
	
	Solving ever larger problems towards practical applications requires, at some point, the use of more than one core to speed up the simulations.
	Furthermore, adaptive refinement in the regions of interest is of huge importance.
	This aims to keep the problem size as small as possible without losing accuracy of the final solution.
	On the other hand, adaptivity usually poses additional difficulties towards the implementation and the linear solvers.
	
	
	\subsection{A quasi-monolithic semi-linear form}
	
	For the solution process, we add both equations in \cref{problem} and define a common semi-linear form $A$.
	Here, we use a linear extrapolation \cite{HeWhWi15} in the first term
	of the $u$-equation in order to deal with the otherwise non-convex
	problem. Thus, we replace $\varphi$ by a time-lagged extrapolation 
		$\tilde\varphi$ yielding $g_\kappa (\tilde\varphi)$ and $\tilde\varphi^2$; see \cite{HeWhWi15}.
	Hence, we obtain 
	\begin{align}
	\begin{split}
	\label{eq:semilinear}
	A(u, \varphi)(w, \psi - \varphi) :=
	&\left(  g_\kappa (\tilde\varphi) \sigma^+(u), e(w) \right) + \left(\sigma^-(u), e(w) \right) + \left(\tilde\varphi^2 p, \mbox{div } w \right) \\
	&+ \frac{1}{2} \left(  \partial_\varphi g_\kappa(\varphi) E_s^+(e(u)), \psi - \varphi \right) + 2 (\varphi \ p \ \mbox{div } u, \psi - \varphi) \\
	&+ G_c \left( -\frac{1}{\varepsilon} ( 1-\varphi, \psi -\varphi ) + \varepsilon \left( \nabla \varphi, \nabla (\psi - \varphi ) \right) \right).
	\end{split}
	\end{align}
	Hence, \cref{problem} reduces to the solution of $A(u, \varphi)(w, \psi - \varphi) \geq 0$.
	In order to deal with the variational inequality, 
	the constraint $\varphi \leq \varphi^{old}$ (prescribed in the set
	$W_{in}$ in \eqref{W_in}) 
	is treated with a primal-dual Active-Set algorithm, which is described in \cref{sec:active_set}.
	
	\begin{remark} 
		In the literature, the phase-field approximation described
		here is also referred to as AT~2 model (named after Ambrosio/Tortorelli) introduced in
			\cite{AmTo90} for the Mumford-Shah problem and the original 
			variational fracture formulation \cite{BoFrMa00}.
			Changing the last line in \cref{problem} leads to another well
			known model referred to as AT~1, e.g., \cite{BoMaMa14}.
		However, in this work we only consider the AT~2 model.
	\end{remark}
	
	
	\subsection{Discretization}
	
	Further on, we will require discrete counterparts of the semi-linear form (\ref{eq:semilinear}) and its derivative.
	For the discretization we use a quadrilateral ($2d$) and hexahedral ($3d$) decomposition of the mesh with $Q_1$ elements for the displacement $u$ and phase-field $\varphi$.
	The implementation was done using the C++ FEM library deal.II \cite{BaHaKa07,AlArBa18}.
	
	With the usual FEM approach, we define the discrete version of \eqref{eq:semilinear} as  $$A_h(\underline{U}_h)_i := A(U_h)(\Phi_i) \quad \forall i = 1, 2, \dots .$$
	Here, $\Phi_i$ denotes the $i$-th shape-function, corresponding to either $u$ or $\varphi$.
	The current FE-function, at which the semi-linear form is evaluated, is given by $U_h := (u_h, \varphi_h)$.
	Its coefficients are given by $\underline{U}_h$, i.e. $U_h = \sum_j \underline{U}_h^j \Phi_j$.
	
	Similarly, we define the discretization of the Jacobian $$G(U_h)_{ij} := A'(U_h)(\Phi_i, \Phi_j).$$
	
	For the sake of simplicity, we omit the underline and use $U$ instead
	of $\underline{U}$ in the remainder of this paper.
	
	\begin{problem}
		\label{discreteproblem}
		Let $V_h \subset V, W_h \subset W$ be standard conforming FE-spaces, 
		and let $p \in L^{\infty}(D)$ be given.
		Find discrete vector-valued displacements and a scalar-valued phase-field variable $\{ u_h, \varphi_h \} := \{ u^n_h, \varphi^n_h \} \in \{u_D + V_h \} \times W_{h,in}$ such that:
		\begin{align*}
		& \left(  g_\kappa (\varphi_h) \sigma^+(u_h), e(w_h) \right) + \left(\sigma^-(u_h), e(w_h) \right) + \left(\varphi_h^2 p, \mbox{div } w_h \right) = 0 \quad & \forall w_h \in V_h, \\
		& & \\
		& \frac{1}{2} \left(  \partial_\varphi g_\kappa(\varphi_h) E_s^+(e(u_h)), \psi_h - \varphi_h \right) + 2 (\varphi_h \ p \ \mbox{div } u_h, \psi_h - \varphi_h) & \\
		& + G_c \left( -\frac{1}{\varepsilon} ( 1-\varphi_h, \psi_h - \varphi_h ) + \varepsilon \left( \nabla \varphi_h, \nabla (\psi_h - \varphi_h ) \right) \right) \geq 0 \quad & \forall \psi_h \in W_h,
		\end{align*}
		at each incremental step $n = 1, 2, 3, \ldots$.
	\end{problem}
	%
	
	
	\section{Solution Approach}
	\label{sec:solution}
	
	In this section, we discuss the question on how to solve \cref{discreteproblem}.
	The remaining key challenge is the variational inequality resulting from 
	the irreversibility condition $\varphi \leq \varphi^{old}$.
	Treating this inequality requires specially tailored solution methods.
	In this work, we chose a variant of the well-known Primal-Dual Active-Set method \cite{ItKu00, HiItKu03}, which we briefly describe in the next section.
	
	The linear systems of equations arising from the Active-Set method will be solved using the GMRES solver \cite{SaSc86} with geometric multigrid preconditioning.
	\Cref{sec:multigrid} describes the details of the multigrid scheme.
	
	\subsection{Primal-Dual Active-Set}
	\label{sec:active_set}
	
	We now turn to the description of the nonlinear solver.
	We deal with two types of nonlinearities: nonlinear equations 
	and nonlinear behavior due to the inequality constraint.
	The goal is to design an algorithm that treats both
	simultaneously\footnote{In this paper, the combined nonlinear
		(semi-smooth) Newton solver
		is called primal-dual active set or, sometimes shorter, as active set (\cref{sec:results}).}.
	The original $2d$ version adapted to the phase-field fracture problem was developed in \cite{HeWhWi15}.
	In \cite{LeWhWi16}, and more recently in \cite{HeWi18}, it has been 
	proven to be computationally reliable also for the $3d$ setting.
	In the following, we omit the subscript $h$ to make the notation more readable.
	
	The general algorithm reads:
	\begin{algorithm}[H]
		\caption{Primal-dual active set}
		\label{alg:active_set}
		Repeat for $k=0, \dots$ until $\mathcal{A}_k = \mathcal{A}_{k-1}$ and $\Vert \widetilde{R}_k \Vert \leq \varepsilon_{as}$:
		\begin{algorithmic}[1]
			\State Assemble residual $R_k = -A(U_k)$
			
			\State Assemble matrix $G_k = A'(U_k)$
			
			\State Compute active set $\mathcal{A}_k = \{i \mid (M^{-1} R_k)_i + c \ (U_k - U^{old})_i > 0 \text{ and } i \sim \varphi \}$
			
			\State Eliminate rows and columns in $\mathcal{A}_k$ from $G_k$ and $R_k$ to obtain $\widetilde{G}_k$ and $\widetilde{R}_k$
			
			\State Solve linear system $\widetilde{G}_k \cdot \delta U_k = \widetilde{R}_k$
			
			\State Update $U_{k+1} = U_k + \delta U_k$

		\end{algorithmic}
	\end{algorithm}
	
	\begin{remark}
		Step 3 requires the computation of $M^{-1} R_k$, where $M$ denotes the usual mass-matrix.
		For practical purposes, we assemble $M$ using a Gauss-Lobatto quadrature rule.
		This leads to a diagonal mass matrix, which is trivial to invert and store.
		Note that the active set acts only on dofs associated to the phase-field $\varphi$ (indicated by $i \sim \varphi$).
		Hence, it is sufficient to perform the computations in steps 3 and 4 only on the phase-field dofs.
		The entries associated to the displacement remain unchanged.
		In our simulations we chose $c = 100$ and $\varepsilon_{as} = 10^{-10}$.
	\end{remark}

	\begin{remark}
		\Cref{alg:active_set} presented above switches between the computation of $\mathcal{A}_k$ (step 3) and a single Newton correction step on the inactive set $\mathcal{A}_k^c$ (steps 5-6).
		We could also modify this method to perform more Newton corrections before updating the active set again.
		However, in our experiments we observed that the bottleneck in
		this method is the convergence of the active set, i.e.,
		$\mathcal{A}_k = \mathcal{A}_{k-1}$ 
		and not the Newton convergence criterion (see also Fig. 14 in \cite{HeWhWi15}).
	\end{remark}
	
	
	\subsection{Geometric Multigrid}
	\label{sec:multigrid}
	
	In this work, we solve the linear equations in Step $5$ in Algorithm
	\ref{alg:active_set} 
	from the Active-Set algorithm using a geometric multigrid approach.
	The main idea of multigrid methods is to solve problems on coarser
	grids and use 
	this information to enhance the solution on the finest grid.
	
	Transfer of vectors from one grid level $L$ to another level $l$ is achieved by the canonical restriction and prolongation operators $\mathcal{R}^L_l$ and $\mathcal{P}^L_l$.
	However, restricting a function to a coarser grid is subject to heavy losses of accuracy, since we cannot represent highly oscillating parts anymore.
	One of the key steps to make multigrid methods work is to "smooth" the vector that needs to be transferred.
	This reduces the oscillating parts and hence, it can be approximated on a coarser grid without loosing too much information.
	For more details regarding multigrid methods we refer the reader to the standard literature, see e.g. \cite{Br93,Ha94}.
	
	Smoothing is handled by the Chebyshev-Jacobi method, which will be described in \cref{sec:smoother}.
	A slight variation of the smoother is used for the solution on the coarsest grid as explained in \cref{sec:coarse}.
	\Cref{sec:precondition} illustrates the use of multigrid methods for systems with block-structure.
	
	\subsubsection{Chebyshev-Jacobi Smoother}
	\label{sec:smoother}
	
	Our final goal is to solve \cref{discreteproblem} (i.e., Step 5 in
		Algorithm \ref{alg:active_set}) in a matrix-free fashion.
	Hence, we do not have the matrix $G$ at our disposal, which narrows the selection of applicable smoothers.
	A frequently used smoother in the context of highly parallel matrix-free computations is the Chebyshev-accelerated Jacobi smoother, which works for positive definite and symmetric problems.
	The important property of the Chebyshev-Jacobi smoother with respect to a matrix-free formulation is, that it only relies on matrix-vector products and an estimate of the largest eigenvalue.
	In particular, an explicit representation of the matrix $G$ is not required.
	
	Chebyshev acceleration is applied on top of a general iterative method with iteration matrix $M$, i.e. $$x_{m+1} = M x_m + b$$ for the solution of $A x = b$.
	The accelerated iteration matrix is then given by $p(M)$, with $p(x)$ denoting the (scaled and shifted) Chebyshev polynomials on the interval $[a, b]$, resulting in $$x_{m+1} = p(M) x_m + b.$$
	For practical purposes, it makes sense to utilize the well-known three-term recurrence relation for Chebyshev polynomials.
	More details can be found in \cite{Va99} and the references provided therein.
	
	The interval $[a,b]$ determines which eigenvalues should be targeted by the acceleration technique.
	For example: if it is used as a solver, we would like to reduce the error with respect to all eigenvalues/vectors.
	Hence, it would be beneficial to have $a := \lambda_{min}, b := \lambda_{max}$.
	However, information about the extremal eigenvalues is (most often) not available.
	Luckily, $\lambda_{max}$ is easy to approximate, e.g. by means of a CG method.
	In our work, the main purpose for the Chebyshev method is to act as a smoother inside a multigrid scheme.
	There, the small eigenvalues are usually handled by the coarse-grid correction.
	Hence, we want to focus the action of the Chebyshev acceleration to some of the largest eigenvalues only, i.e. $[c \ \lambda_{max}, \lambda_{max}]$, with $0 < c < 1$.
	
	In our application, we apply this acceleration technique to the Jacobi method, whose iteration matrix is given by $M = (I - D^{-1} A)$.
	The matrix diagonal $D$ required here can be easily computed also in
	the matrix-free context, 
	as discussed in \cref{sec:matrixfree}. 
	A comparison of different smoothers regarding their performance 
	in HPC (high performance computing) applications is
	demonstrated in \cite{AdBrHu03,BaFaKo11}.

	Within our numerical experiments presented in Section \ref{sec:results}, we approximate $\lambda_{max}$ using $10$ CG iterations.
	A safety factor of $1.2$ is included, to account for possible underestimation of $\lambda_{max}$ by the CG method.
	This is required, since the total number of iterations increases a lot, if the smoothing range is smaller than $\lambda_{max}$.
	On the other hand, overestimating $\lambda_{max}$ does not significantly alter the number of required iterations.
	This behavior is also observed in \cite{AdBrHu03}.
	As denoted before, we chose the lower eigenvalue bound as a fraction of the upper bound.
	Summarizing, we define the smoothing range as $[1.2 \ \tilde\lambda_{max} / 5, 1.2 \ \tilde\lambda_{max}]$, where $\tilde\lambda_{max}$ is the approximation obtained via the CG method.
	Furthermore, we restrict ourselves to Chebyshev polynomials of degree $4$.
	
	\subsubsection{Coarse-Grid Solver}
	\label{sec:coarse}
	
	On the coarsest multigrid level, we use the Chebyshev method again but now for solving instead of smoothing.
	We adjust the "smoothing" range to $[0.9 \ \tilde{\lambda}_{min}, 1.2 \ \tilde\lambda_{max}]$, i.e. treating the whole (approximate) spectrum.
	In our experiments, using more expensive solvers (e.g. LU-factorization) did not yield significant improvements.
	
	\subsection{Preconditioning}
	\label{sec:precondition}
	
	We recall that at Step 5 in Algorithm \ref{alg:active_set} we
	solve
	\[
	\widetilde{G} \cdot \delta U = \widetilde{R}
	\]
	where $\widetilde{G}$ has the block structure
	\[
	\widetilde{G} = 
	\begin{pmatrix}
	G_{uu} & G_{u\varphi}\\
	G_{\varphi u} & G_{\varphi\varphi}
	\end{pmatrix}
	\]
	The block $G_{u\varphi}$ is zero because we use the time-lagged
	extrapolation in the displacement equation and therefore the
	linearization with respect to $\varphi$ is zero.
	
	The preconditioned system reads:
	\[
	P^{-1} \widetilde{G} \cdot \delta U = P^{-1} \widetilde{R}
	\]
	with a preconditioner $P^{-1}$ approximating the inverse of $\widetilde{G}$.

	With these preparations, the geometric multigrid can be applied in several ways.
	One possibility is to utilize the block-structure in matrix $\widetilde{G}$ and apply a multigrid scheme to each of the diagonal blocks.
	This yields the block-diagonal preconditioner $$P_{diag}^{-1} := \begin{bmatrix} MG(\widetilde{U}) & 0 \\ 0 & MG(\widetilde{P}) \end{bmatrix},$$ with Chebyshev-Jacobi smoothers $\mathcal{S}$ used inside each of the multigrid methods $MG$.
	
	A different approach is to apply the multigrid scheme to the entire linear system $\widetilde{G}$, yielding $$P_{full}^{-1} := MG(\widetilde{G}).$$
	Within this approach, the block-structure is utilized inside the smoother in a similar manner as before:
	$$\mathcal{S}_{full} := \begin{bmatrix} \mathcal{S}(\widetilde{U}) & 0 \\ 0 & \mathcal{S}(\widetilde{P}) \end{bmatrix},$$ with Cheybshev smoothers for the blocks $U$ and $P$.
	One could extend the diagonal smoother $\mathcal{S}_{full}$ by adding the off-diagonal block.
	However, we did not observe a significant decrease in the number of required iterations.
	
	
	\section{Matrix-Free}
	\label{sec:matrixfree}
	
	Due to ever increasing demand to solve larger and larger problems, memory limitations tend to become more of an issue.
	To avoid these limitations, we aim to implement the geometric multigrid solver in a matrix-free fashion.
	This allows us to handle large problems without the need for enormous amounts of memory.
	In this section, we describe some basic aspects and algorithmic details regarding the matrix-free realization of the geometric multigrid solver.
	The main components are the matrix-free vector multiplication on the finest level and the corresponding part on the coarser levels.
	We split the discussion into these two parts, where the former contains common aspects of the matrix-free framework, whereas the latter contains some peculiarities regarding the operators on the coarser levels.
	
	\subsection{Matrix-Vector Multiplication in the Matrix-Free Context}
	
	In this part, we briefly illustrate the main principles behind matrix-free techniques.
	More details on matrix-free and its integration into deal.II are found in \cite{KrKo12}, and \cite{KrKo17} for discontinuous elements.
	
	Within the matrix-free context, we want to evaluate $v = G \cdot u$ without actually assembling the FE-matrix $G := G(U_h)$.
	This product may be rewritten in an element-wise fashion $$G \cdot u = \sum_{k=1}^{n_e} C^T P_k^T G_k (P_k C u).$$
	Here, $n_e$ denotes the number of elements, $C$ defines possible constraints on dofs, $P_k$ is the element-wise global-to-local mapping and $G_k$ are the local stiffness matrices.

	We can further rewrite the local contributions by transforming it to the reference element.
	For simplicity, we illustrate this for the Laplace problem $(\kappa(x) \nabla u , \nabla v )$.
	Then we have $$G_k = B^T_R J^{-1}_k D_k J_k^{-T} B_R,$$ with $J_k$ denoting the reference element transformation, $B_R$ the gradient matrix on the reference element and $D_k$ a diagonal matrix representing the coefficient $\kappa$.
	The precise definition of all quantities is given in \cite{KrKo12}.
	An important further optimization consists of utilizing the tensor-product structure for computation of $B_R$ on the reference element $[0,1]^d$.
	This changes the computational complexity from $\mathcal{O}(p^{2d})$ to $\mathcal{O}(p^{d+1} d^2)$ per element.
	
	\begin{remark}
		Similar derivations can be carried out for different bilinear forms.
		Luckily, this is already included in deal.II for operators of the form $(\cdot, \Phi)$, $(\cdot, \nabla \Phi)$ and $(\cdot, \Delta \Phi)$, i.e. arbitrary terms tested with the test-function $\Phi$, its gradient or its Laplacian.
		This is more than sufficient to cover most partial differential equations.
		Furthermore, deal.II implements vectorization with respect to the elements, i.e. computations are carried out on multiple elements at the same time, if run on a reasonably modern CPU \cite{KrKo12,KrKo17}.
	\end{remark}
	
	\begin{remark}
		In a similar manner, one can explicitly compute the vector of diagonal entries of the matrix $A$.
		This is required for the Jacobi-smoother in the multigrid method.
		Furthermore, the matrix-free framework can also be used to compute the residual $R_k$ required within the nonlinear solver. 
	\end{remark}
	
	\subsection{Coarse-Level Operators}
	
	The operators on the coarse levels are implemented in a matrix-free fashion as well.
	Most parts of the implementation matches those on the fine level $L$, however, there are some important differences.
	First of all, the active set, linearization point and old solutions naturally live on the finest grid.
	In order to apply the matrix-free operator, we need to transfer this information to the coarser levels $l$ as well.
	By $\mathcal{R}_l$ we will denote the restriction operator mapping from the finest level to level $l$.
	
	For the active set, our approach starts by building a vector $a$ representing the active set $\mathcal{I}_\mathcal{A}$ by 
	$$ a := (a^i) := \begin{cases}  1 &i \in \mathcal{A} \\ 0 &\text{else} \end{cases} $$
	for all dof indices $i$.
	The active sets on the coarser levels $\mathcal{A}_l$ are then computed via the restrictions $a_l := \mathcal{R}_l a$ by $\mathcal{A}_l := \{ i : a_l^i = 1 \}.$
	Note that the entries $a_l^i \in \{0, 1\}$, since we are using $Q_1$ elements on nested grids, i.e. each vertex (and hence dof) on level $l$ also exists on all finer levels.
	The big advantage of this method is that it only uses the restriction operators $\mathcal{R}_l$, which are readily available in the multigrid context.
	
	The natural way to evaluate the coarse level operators would include information about the linearization point and old solutions from the finest level, i.e. we would need to evaluate $G_l(U) \cdot v_l$ analogously to $G(U) \cdot v$ on the finest level as before.
	However, this is very expensive, since for each quadrature point on the coarse levels we have to
	\begin{enumerate}
		\item find the fine-grid element this point belongs to,
		\item evaluate all fine-level basis functions (or its gradients) on this element, and
		\item compute the value/gradient using the current linearization point on the fine level.
	\end{enumerate}
	One could precompute parts of this information in order to save computational time, however, this would  defy the whole idea of matrix-free, as it introduces lots of memory overhead.
	Furthermore, these values would have to be updated again every time the linearization point changes.
	
	To overcome this, we approximate the evaluation of $G_l(U) \cdot v_l$ by $G_l(U_l) \cdot v_l$, with $U_l := \mathcal{R}_l U$ being the restriction to the coarser levels.
	Thus, evaluating $G_l$ now depends only on information from level $l$.
	This is way cheaper, since step (1) is no longer required and step (2) is required anyway for the matrix-free multiplication (but now on the coarse elements), whereas step (3) is cheap in any case.
	The only additional cost comes from the computation and storage of $U_l = \mathcal{R}_l U$ for all levels.
	Again, this computation needs to be redone every time $U$ changes.
	
	\begin{remark}
		For the sake of simplicity, we omitted the dependence of $G$ on the old solutions $U^{n-1}, U^{n-2}$ (required for the extrapolation).
		These can be treated in the very same way as $U := U^n$.
	\end{remark}
	
	
	\section{Numerical Results}
	\label{sec:results}
	
	We proceed to show the capabilities of the presented multigrid method applied to some test cases frequently found in the literature.
	In all numerical experiments, we mainly vary the level of refinement and consider fixed values for (most) parameters.
	The linear solver consists of an outer GMRES solver, preconditioned by the previously presented monolithic multigrid method.
	The stopping criterion for the Active-Set method is determined by an absolute tolerance $\varepsilon_{as} := 10^{-10}$.
	The linear solver stops once the absolute residual is less than $10^{-10}$ or a reduction of $10^{-4}$ compared to the initial (linear) residual is reached, i.e. $\Vert r_k \Vert \leq \mbox{max}\{ 10^{-12}, 10^{-4} \Vert r_0 \Vert \}$.
	
	Quantities of interest include loads and energies.
	The crack energy is given by $$\int_{\Omega} \frac{G_c}{2} \left( \frac{1}{\varepsilon} (1 - \varphi)^2 + \varepsilon (\nabla \varphi, \nabla \varphi ) \right) dx.$$
	For the bulk energy we have $$\int_{\Omega} \left( g_\kappa(\varphi) E_s^+(u) + E_s^-(u) + \varphi^2 p \ \mbox{div } u \right) dx.$$
	The loading (in $y$-direction) is evaluated on parts of the boundary using $$L_y = \int_{\Gamma_L} (\sigma \cdot n)_y ds.$$
	
	In the upcoming sections, we present several test scenarios.
	We start with a setting showing multiple fractures joining and branching, taken from \cite{HeWhWi15,HeWi18}.
	We proceed with a L-shaped panel test, which illustrates crack initiation without predefined fractures.
	This test is frequently found in the literature, see e.g. \cite{AmGeDe15,BeMoCh12,FeHo06,UnEcKo07,Wi17,MeDu07} and was first reported in \cite{Wi01}.
	We conclude with a simple $3d$ variation of the L-shaped test to show the applicability of our method also in the $3d$ setting.
	
	\newcommand{\summary}[1]{
		\pgfplotstabletypeset[
		columns={lvl, dofs, cells, h, eps},
		%
		columns/lvl/.style={column name=$l$},
		columns/dofs/.style={sci, zerofill, precision = 1},
		columns/h/.style={column name=$h$, sci, precision = 1},
		columns/eps/.style={column name=$\varepsilon$, sci, precision = 1},
		columns/cells/.style={column name={elements}, sci, zerofill, precision = 1},
		every head row/.style={before row=\toprule,after row=\midrule},
		every last row/.style={after row=\bottomrule},
		]{#1}
	}
	
	\clearpage
	\subsection{Multiple Fractures}
	
	\subsubsection{Problem Description}
	
\begin{figure}[ht]
		\centering
		\begin{tikzpicture}[
		acteur/.style={
			circle,
			draw=black,
			fill=black,
			inner sep=1pt,
			minimum size=0cm
		}
		]
		
		\tikzmath{\h = 0.05; \d = 0.3;}
		
		\node[acteur, label={[below]\tiny $(0, 0)$}] (A) at (0,0) {};
		\node[acteur, label={[above]\tiny $(0, 4)$}] (B) at (0,4) {};
		\node[acteur, label={[above]\tiny $(4, 4)$}] (C) at (4,4) {};
		\node[acteur, label={[below]\tiny $(4, 0)$}] (D) at (4,0) {};
		
		\draw (A) -- (B) -- (C) -- (D) -- (A);
		
		\coordinate (V1) at (2.5 - \h , 0.8);
		\coordinate (V2) at (2.5 + \h , 0.8);
		\coordinate (V3) at (2.5 + \h , 1.5);
		\coordinate (V4) at (2.5 - \h , 1.5);
		
		\draw (V1) -- (V2) -- (V3) -- (V4) -- (V1);
		\draw[|-|] ($ (V1) + (\d, 0) $) -- ($ (V4) + (\d, 0) $);
		\node[fill=none, label={[right, anchor=north west]\tiny{$0.7$}}] at ($ (V1) ! 1/2 ! (V4) + (\d, 0.1)$) {};	
		\draw[|-] ($ (V3) + (0, \d/2)$) -- ($ (V3)  + (0.15, \d/2)$);
		\draw[|-] ($ (V4) + (0, \d/2)$) -- ($ (V4)  + (-0.15, \d/2)$);
		\node[fill=none, label={[above, yshift=-1mm]\tiny $\varepsilon$}] at ($ (V3) ! 1/2 ! (V4) + (0, \d/2)$) {};
		
		\coordinate (H1) at (0.5 , 3.0 + \h);
		\coordinate (H2) at (0.5 , 3.0 - \h);
		\coordinate (H3) at (1.5 , 3.0 - \h);
		\coordinate (H4) at (1.5 , 3.0 + \h);
		
		\draw [] (H1) -- (H2) -- (H3) -- (H4) -- (H1);
		\draw[|-|] ($ (H2) + (0, \d) $) -- ($ (H3) + (0, \d) $);
		\node[fill=none, label={[above, yshift=-1mm]\tiny{$1.0$}}] at ($ (H2) ! 1/2 ! (H3) + (0, \d)$) {};
		\draw[|-] ($ (H3) + (\d/2, 0)$) -- ($ (H3)  + (\d/2, -0.15)$);
		\draw[|-] ($ (H4) + (\d/2, 0)$) -- ($ (H4)  + (\d/2, +0.15)$);
		\node[fill=none, label={[right, yshift=-1.5mm]\tiny $\varepsilon$}] at ($ (H3) ! 1/2 ! (H4) + (\d/2, 0)$) {};
		\end{tikzpicture}
		
		\caption{Multiple Fractures.
			Geometry and position of initial cracks.}
		
		\label{multiple:geometry}
	\end{figure}
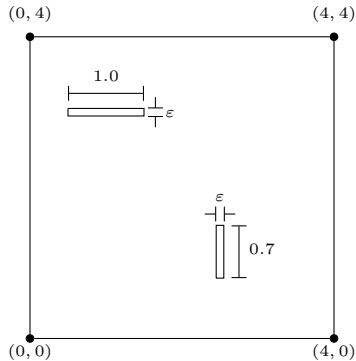

	This problem consists of two prescribed fractures inside the domain $\Omega = ( 0\mbox{ m} , 4\mbox{ m} ) ^2$ as shown in \cref{multiple:geometry}. 
	Further propagation of the fractures is driven by an increasing pressure inside the cracks, in particular $p(t) = 10^3 t \mbox{ Pa}$.
	For the computation, a constant time increment of $dt := 10^{-2} s$ has been used.
	The initial cracks are given at $ \{ 2.5 - h / 2 \leq x \leq 2.5 + h / 2 \text{ and } 0.8 \leq y \leq 1.5 \} $ (vertical fracture) and  $ \{  0.5 \leq x \leq 1.5   \text{ and } 3 - h / 2 \leq y \leq 3 + h / 2 \} $ (horizontal fracture), with the element diameter $h$.
	Initial conditions are given by homogeneous Dirichlet conditions for $u$.
	The phase-field regularization parameter is chosen as $\varepsilon :=
	h$.
	The full list of parameters is summarized in \cref{multiple:summary,multiple:parameters}.
	
	\begin{table}[ht]
		\renewcommand{\arraystretch}{1.3}
		\centering
		\summary{data/summary_multiple_refinement_heps/summary.csv}

		\caption{Multiple Fractures.
			Values and parameters for the different tests.}
		
		\label{multiple:summary}
	\end{table} 
	
	\begin{table}[ht]
		\renewcommand{\arraystretch}{1.3}
		\centering
		\begin{tabular}{c|c|l}
			Variable & Value & Unit \\ \hline
			$\nu$ & $0.2$ & $[1]$ \\
			$E$ & $10^4$ & $[Pa]$\\
			$\mu$ & $ 4166.6 $ & $[N/m^2]$ \\
			$\lambda$ & $ 2777.7 $ & $[N/m^2]$ \\
			$\alpha$ & $0.0$ & $[1]$ \\
			$p(t)$ & $10^3 t$& $[N/m^2]$ \\
			$G_c$ & $1.0$ & $[N / m]$ \\
			$\kappa$ & $10^{-10}$ & $[1]$ \\
			$\varepsilon$ & $0.044$ & $[m]$ \\
			$dt  \ [s]$ & $10^{-2}$ & $[s]$ \\
		\end{tabular}
		
		\caption{Multiple Fractures.
			Material parameters and configuration. Lam\'e parameters $\mu, \lambda$ are derived from the Poisson ratio $\nu$ and elastic modulus $E$.}
		
		\label{multiple:parameters}
	\end{table}
	
	\subsubsection{Random Variations}
	
	Following the original experiments presented in \cite{HeWhWi15}, we consider a randomly perturbed elastic modulus for the multiple fracture problem.
	The elastic modulus is varied between $0.1$ and $1$-times its original value.
	The degree of perturbation is given by a "smooth" random field $r(x,y)$, illustrated in \cref{multiple:pictures}.
	
	Computation of $r(x,y)$ was done using the Simplex algorithm provided by the C++ library FastNoise \cite{Au18}.
	The precise settings to reproduce the created random field are (without explanation): seed $= 2$, frequency $= 0.5$, simplex noise and linear interpolation.
	Source code and an explanation of the different settings can be found in the online repository\footnote{https://github.com/Auburns/FastNoise, as of Dec. 1st, 2018}.
	
	This provides a way of deterministic evaluation of the random field at given coordinates, i.e. equal coordinates always result in the same value.
	Therefore, this information can be used consistently on all multigrid levels and possibly multiple cores without additional effort.
	
	\begin{remark}
		The main purpose of these random variations is to show that our solver is capable of handling perturbations.
		We do not aim to provide a realistic representation of inhomogeneous material coefficients.
	\end{remark}
	
	\subsubsection{Numerical Results}
	
\begin{figure}[ht]
		\includegraphics{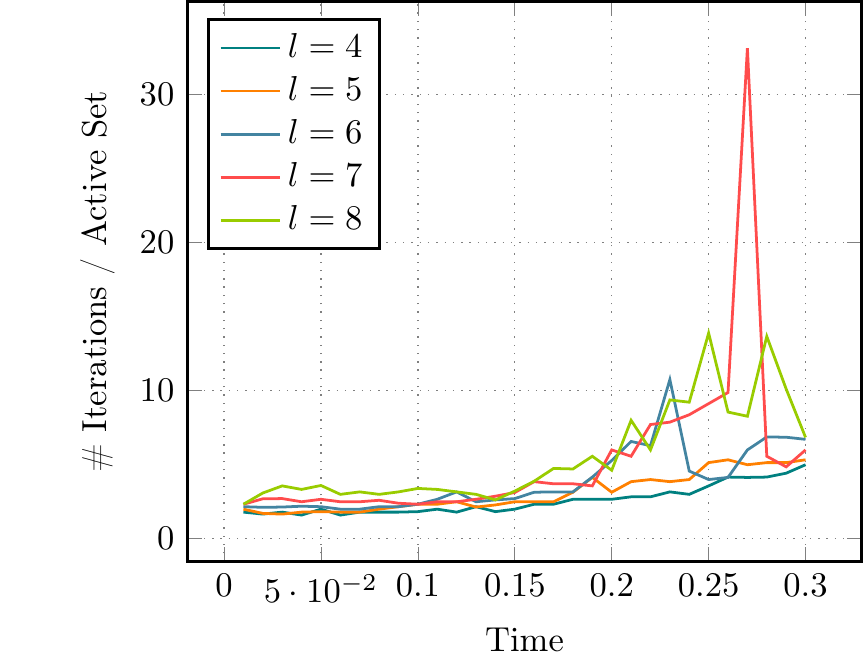}
		\includegraphics{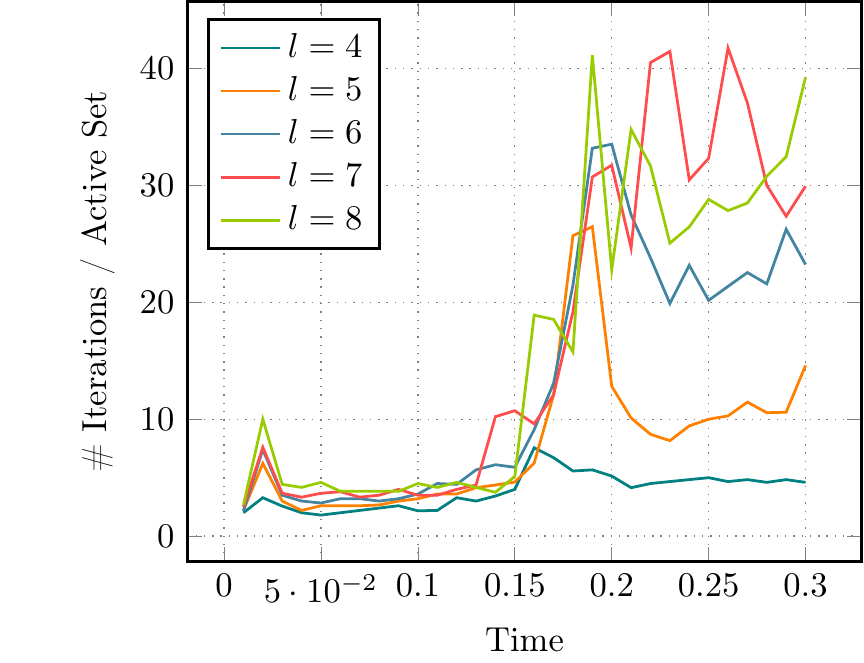}
		
		\caption{Multiple fractures. 
			Iterations of the linear solver per active set step over time for different refinement levels $l$ with $\varepsilon = h$. 
			Left: no stress splitting, right: Miehe-type stress splitting.}
		
		\label{multiple:iterations}
	\end{figure}

	In \cref{multiple:iterations} (left), the number of iterations of the linear solver (GMRES) per Active-Set step is shown for different refinement levels.
	We observe an overall increase in the number of iterations as the fracture grows, with larger spikes when the crack touches the boundary.
	Up to this moment, the effect of $h$-refinement is negligible, showing the robustness of our method.
	Afterwards, the behavior is more irregular, with the largest spike happening for $l=7$.
	
	In the case of Miehe-splitting (\cref{multiple:iterations}, right), we require more iterations of the linear solver.
	Furthermore, the irregular iteration counts become more prominent once the fracture gets close to the boundary ($\sim T = 0.2s$).
	Similar observations are obtained for the Amor-type splitting.
	
\begin{figure}[ht]
		\includegraphics{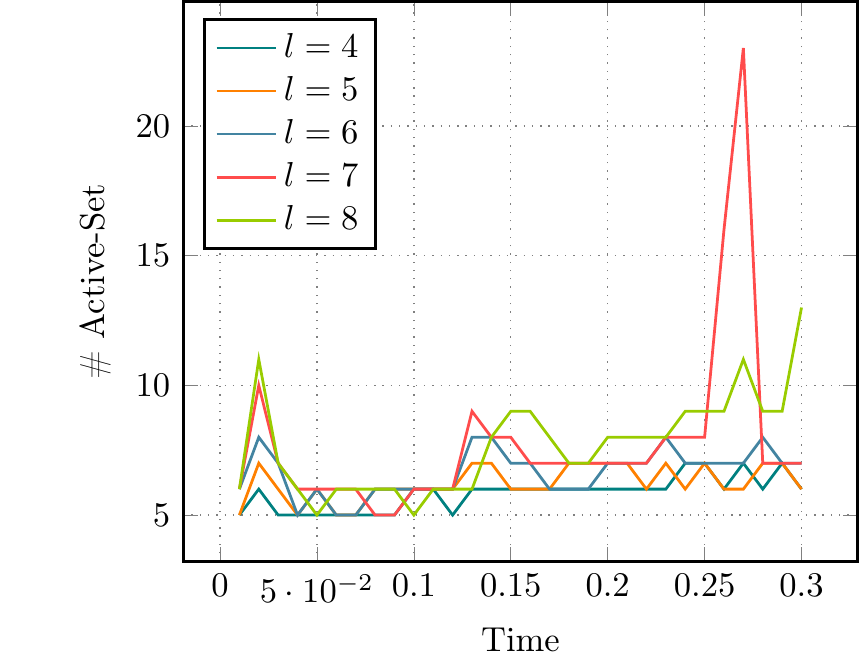}
		\includegraphics{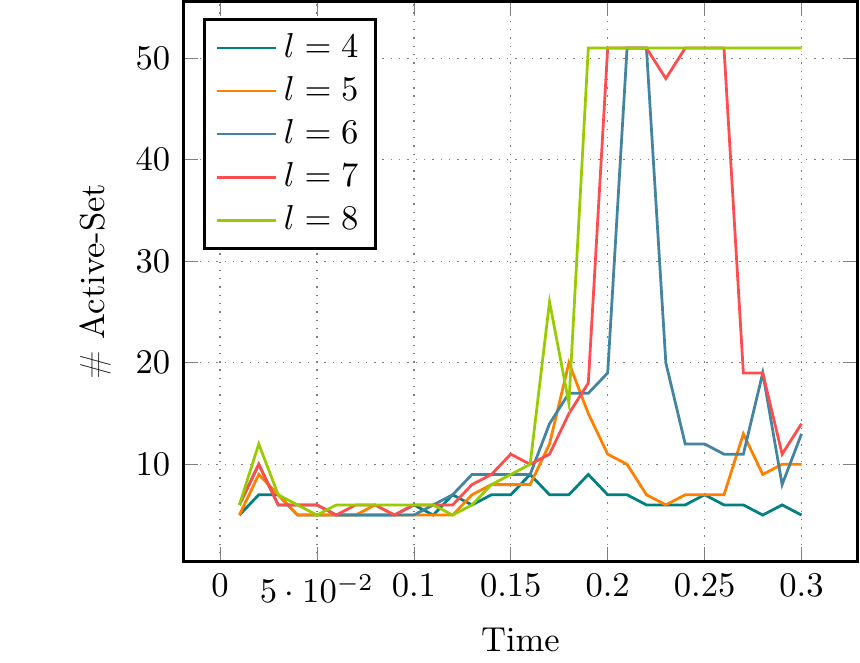}
		
		\caption{Multiple fractures (Example 1).
			Number of Active-Set steps over time for different
			refinement levels $l$ with $\varepsilon = h$. 
			Left: no splitting. Right: Miehe-splitting.
		}
		
		\label{multiple:activeset}
	\end{figure}

	The number of Active-Set iterations behaves quite similar. \Cref{multiple:activeset} shows these numbers over time for different refinement levels. Again, in the case without splitting, the iteration counts remain quite stable, with some peaks towards the end of the simulation. In the Miehe-case, these peaks become more severe.
	
\begin{figure}[ht]
		\centering
		\includegraphics[width=0.25\textwidth]{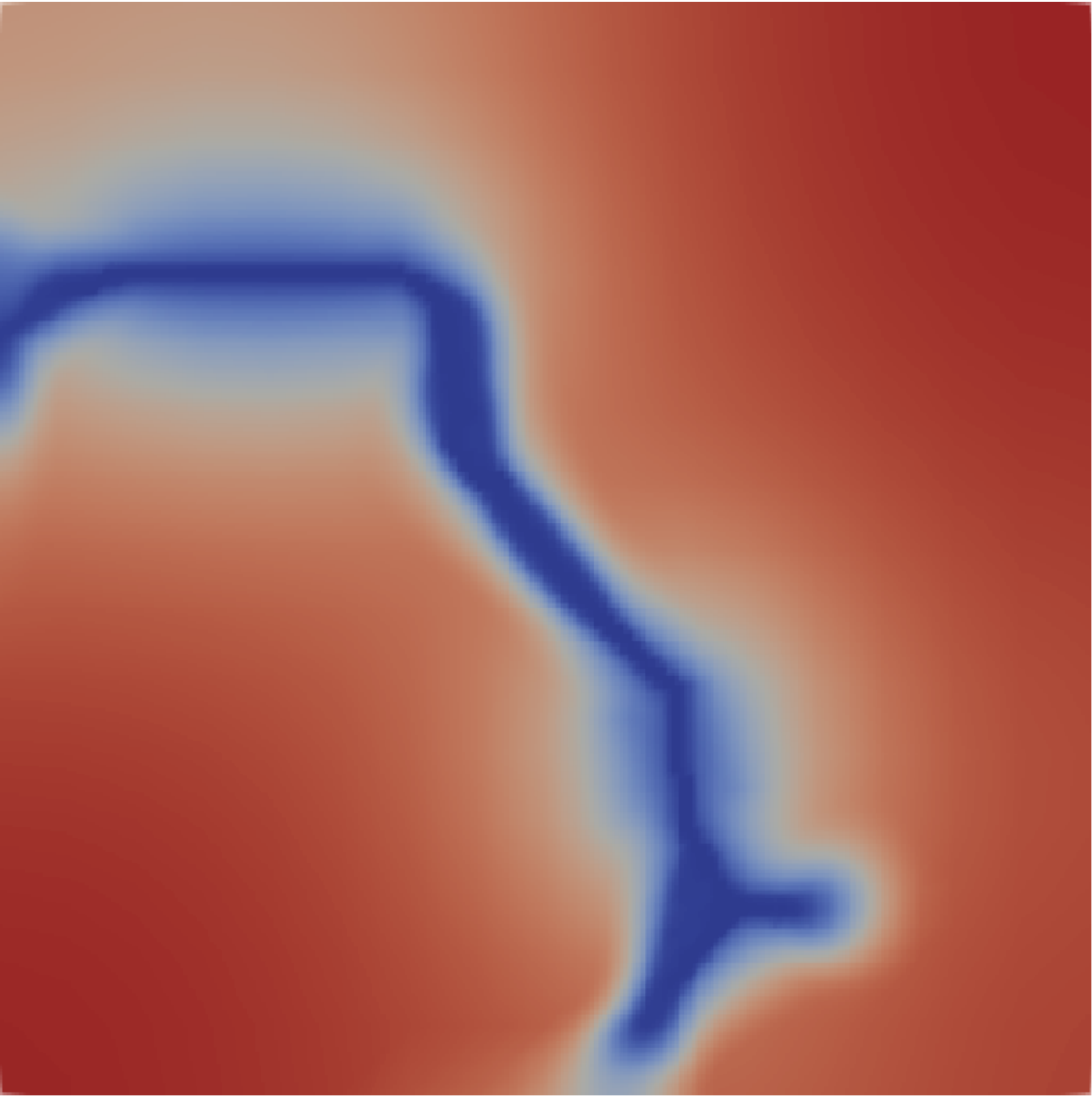}
		\hspace{10mm}
		\includegraphics[width=0.25\textwidth]{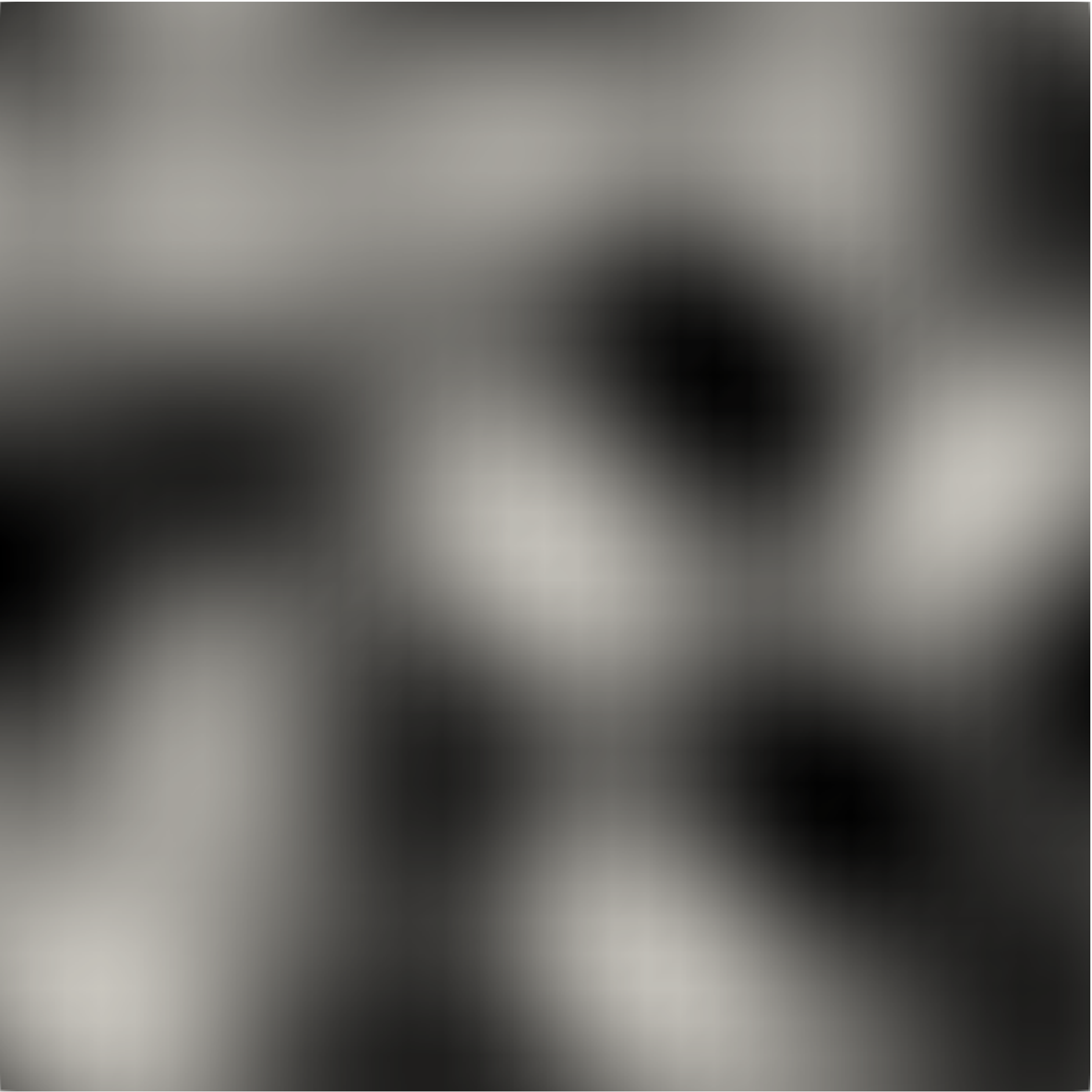}
		
		\caption{Multiple fractures. 
			Fracture pattern using random variations in the elastic modulus ranging from $0.1$ (white) to $1$ (black) times the original value as shown in the right figure.}
		
		\label{multiple:pictures}
	\end{figure}

	\Cref{multiple:pictures} shows the resulting fracture  %
	using randomly perturbed values for the elastic modulus.
	The random variations range from $0.1$ to $1$-times the original value, indicated by the image on the right.
	There, black denotes large values of elastic modulus, whereas white denotes small values.
	It can be observed that the resulting fracture avoids the high elastic modulus regions.

\begin{figure}[ht]
		\centering
		\includegraphics{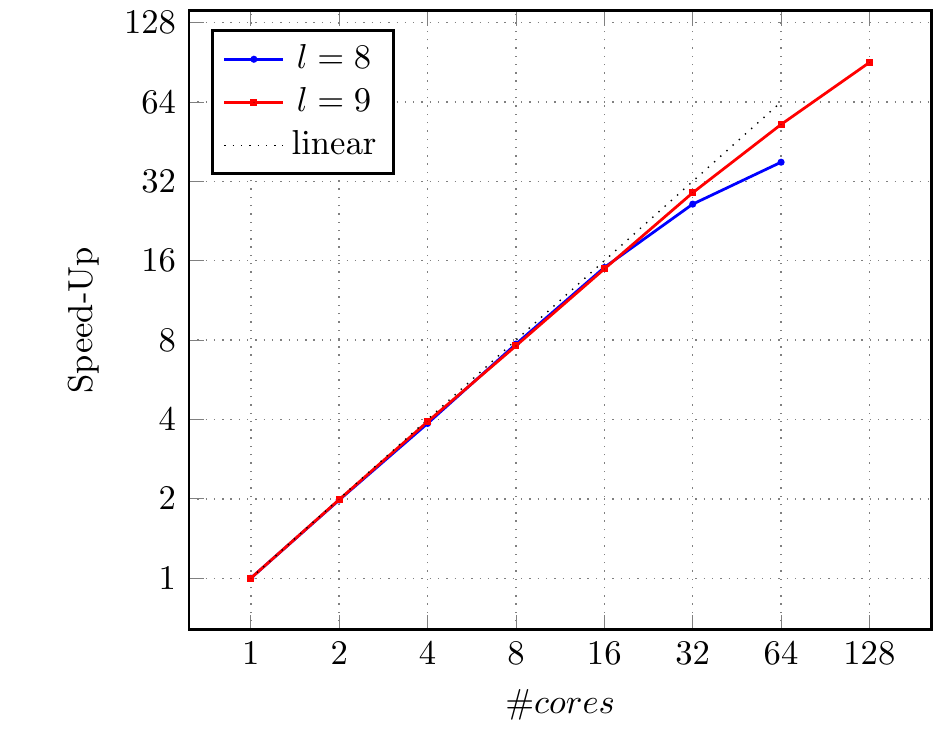}
		
		\caption{Multiple fractures.
			Speedup using distributed parallelization depending on the number of cores for a full simulation with $\varepsilon = h$ and roughly $0.8 m$ and $3.2 m$ dofs.
			Remarkably, the total number of iterations (Active-Set, linear) matches perfectly for all numbers of cores.
		}
		
		\label{multiple:scaling}
	\end{figure}

	Matrix-free implementations are particularly suitable for parallelization.
	A first result into this direction is shown in \cref{multiple:scaling}.
	There, the speed-up for a full simulation of the multiple fracture test with $0.8 m$ and $3.2 m$ dofs is shown.
	The parallel (strong) scaling using distributed parallelization is close to perfect down to a local problem size of roughly $25 k$ dofs.
	When using more cores, the performance drops due to the increased overhead required for communication for smaller local problems.
	This allows us to reduce the computational time for the whole simulation from $400$ minutes down to $10$ minutes on $64$ cores for $l = 7$.
	In the case of $l = 8$, the simulation time drops from $35$ hours down to $23$ minutes using $128$ cores.
	
	\clearpage
	\subsection{L-Shaped Panel}
	
	\subsubsection{Problem Description}
	
	In this section, we consider the L-shaped panel test, which illustrates crack initiation and propagation without any predefined fractures.
	
\begin{figure}[ht]
		\centering
		\begin{tikzpicture}
		[
		acteur/.style={
			circle,
			draw=black,
			fill=black,
			inner sep=1pt,
			minimum size=0cm
		}
		]
		\node[acteur, label={[above, xshift=4mm]\tiny $(0, 0)$}] (A) at (0,0) {};
		\node[acteur, label={[above, xshift=-5mm]\tiny $(250, 0)$}] (B) at (2.5, 0) {};
		\node[acteur, label=above:{\tiny $(250, 250)$}] (C) at (2.5, 2.5) {};
		\node[acteur] (D) at (5, 2.5) {};
		\node[acteur, label=above:{\tiny $(500, 500)$}] (E) at (5, 5) {};
		\node[acteur, label=above:{\tiny $(0, 500)$}] (F) at (0, 5) {};
		\node[acteur, label={[above]\tiny $(470, 250)$}] (U) at (4.4, 2.5) {};
		
		\draw (A) -- (B) -- (C) -- (D) -- (E) -- (F) -- (A);
		\draw[line width=1pt] (U.center) -- (D.center);
		\draw[line width=1pt] (A.center) -- (B.center);
		
		\coordinate (M1) at ($ (U) ! 1/3 ! (D) $);
		\coordinate (M2) at ($ (U) ! 2/3 ! (D) $);
		\draw[<->] (M1) -- ($ (M1) + (0, -0.5) $);
		\draw[<->] (M2) -- ($ (M2) + (0, -0.5) $);
		\node[fill=none, label=below:{\tiny $u_y(t)$}] () at ($ (U) ! 1/2 ! (D) - (0, 0.3) $) {};
		
		\tikzmath{\n = 10;}
		\foreach \i in {0,...,\n}
		\draw ($ (A) ! \i / \n ! (B) $) -- ($ (A) ! \i / \n ! (B) - (0.15, 0.3) $);
		
		\node[fill=none, label=above:{\tiny $\Gamma_{up}$}] () at ($ (E) ! 1/2 ! (F) $) {};
		
		\node[inner sep=0pt, anchor=north west, xshift=1.5cm] (pic) at (E)  {\includegraphics{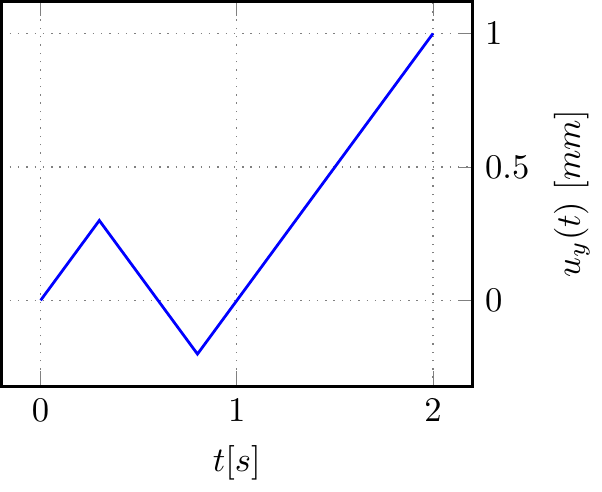}};
		\end{tikzpicture}
		
		\caption{L-shaped panel test. 
			Geometry and cyclic displacement boundary conditions for the L-shaped panel test.}
		
		\label{lshape:configuration}
	\end{figure}
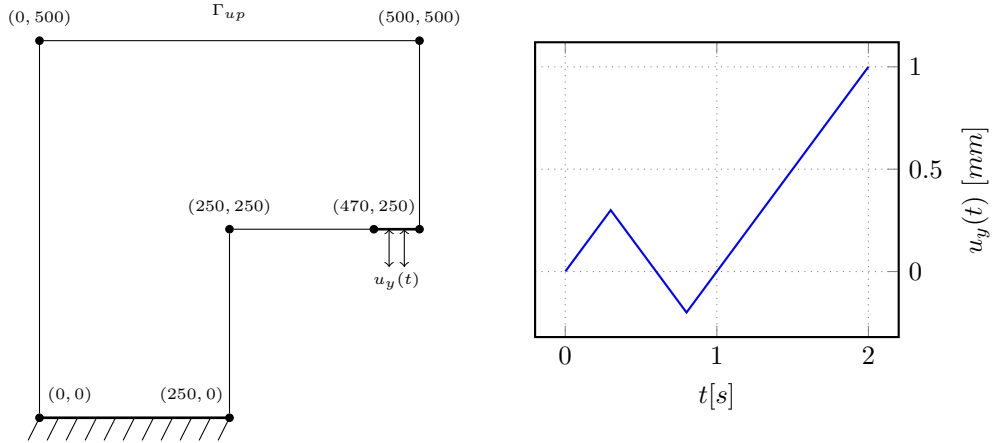

\begin{figure}[ht]
		\centering
		\includegraphics[width=0.2\textwidth]{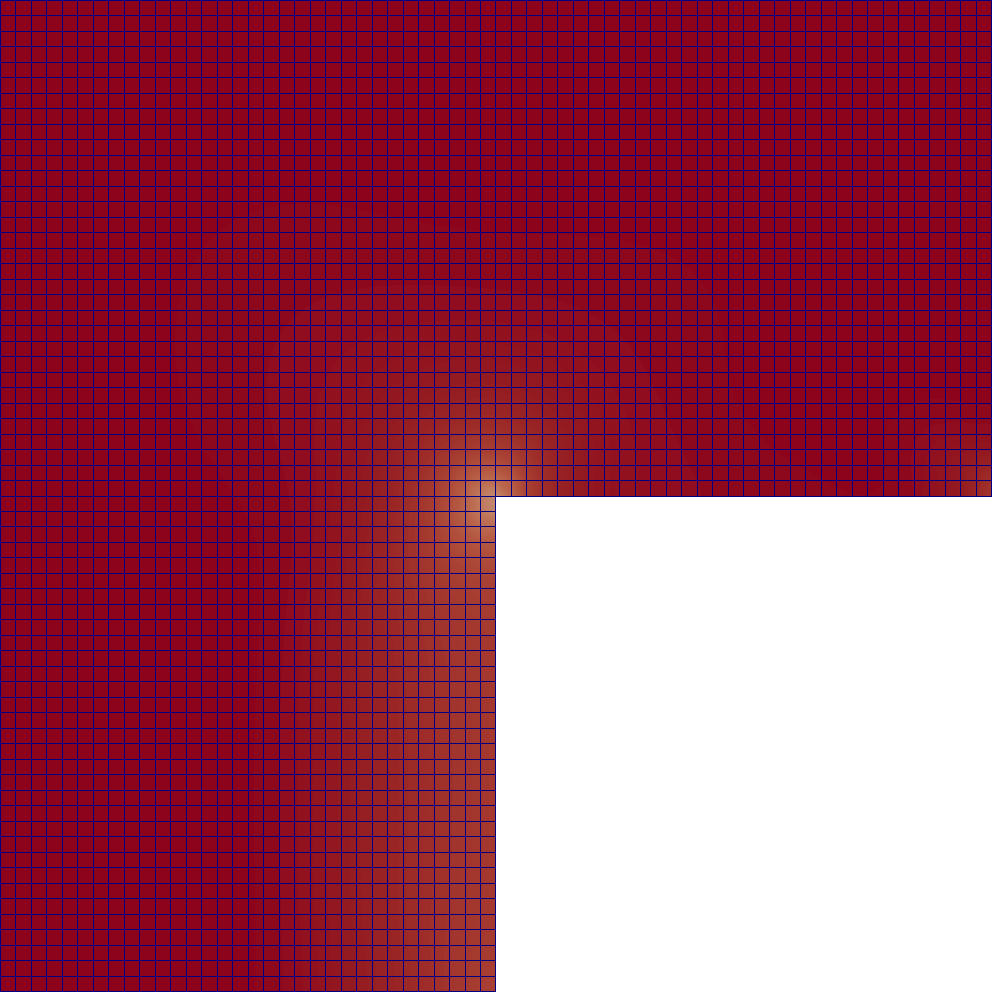}
		\includegraphics[width=0.2\textwidth]{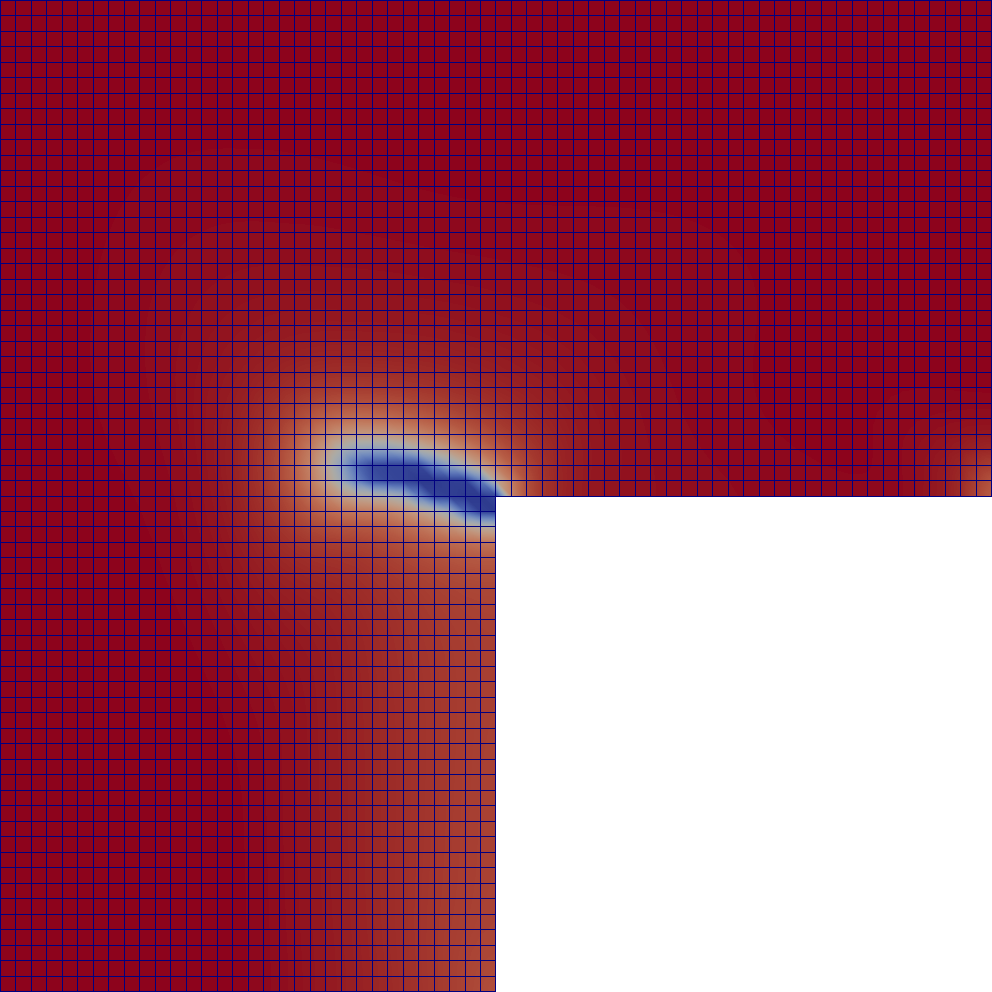}
		\includegraphics[width=0.2\textwidth]{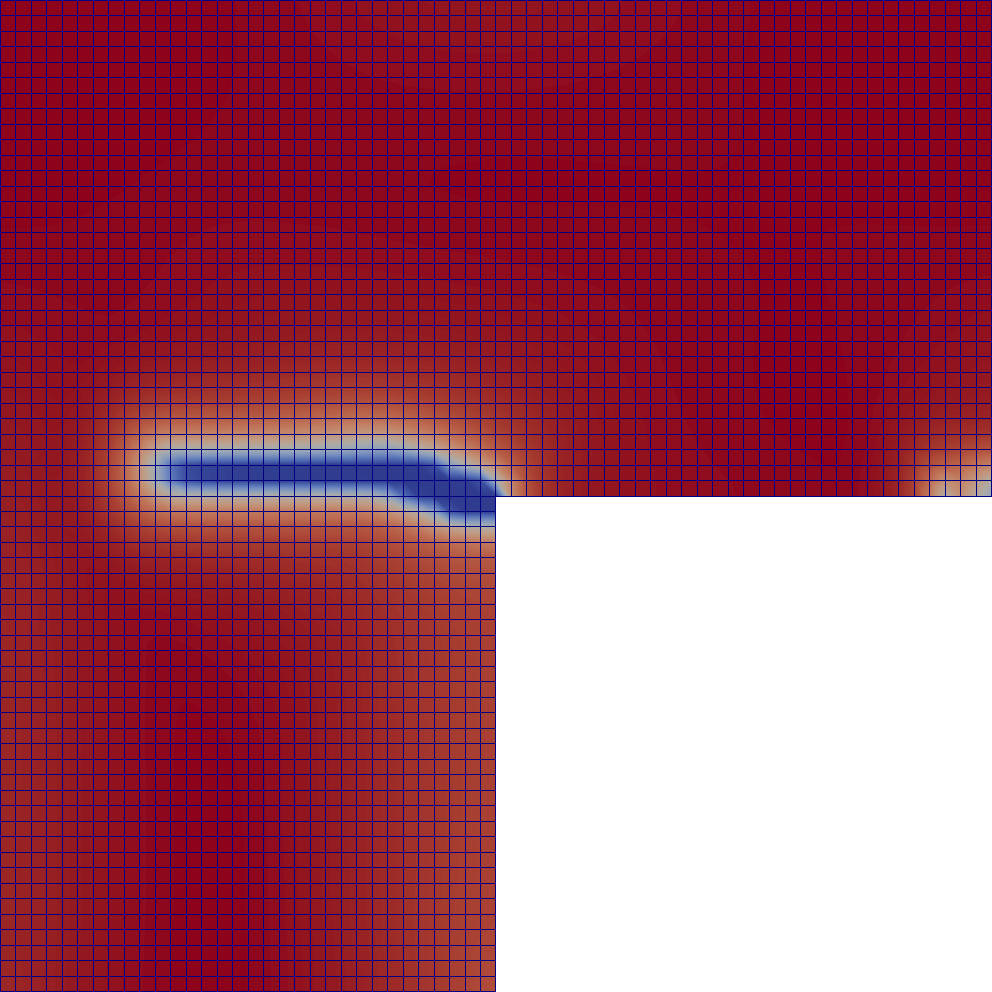}
		\includegraphics[width=0.2\textwidth]{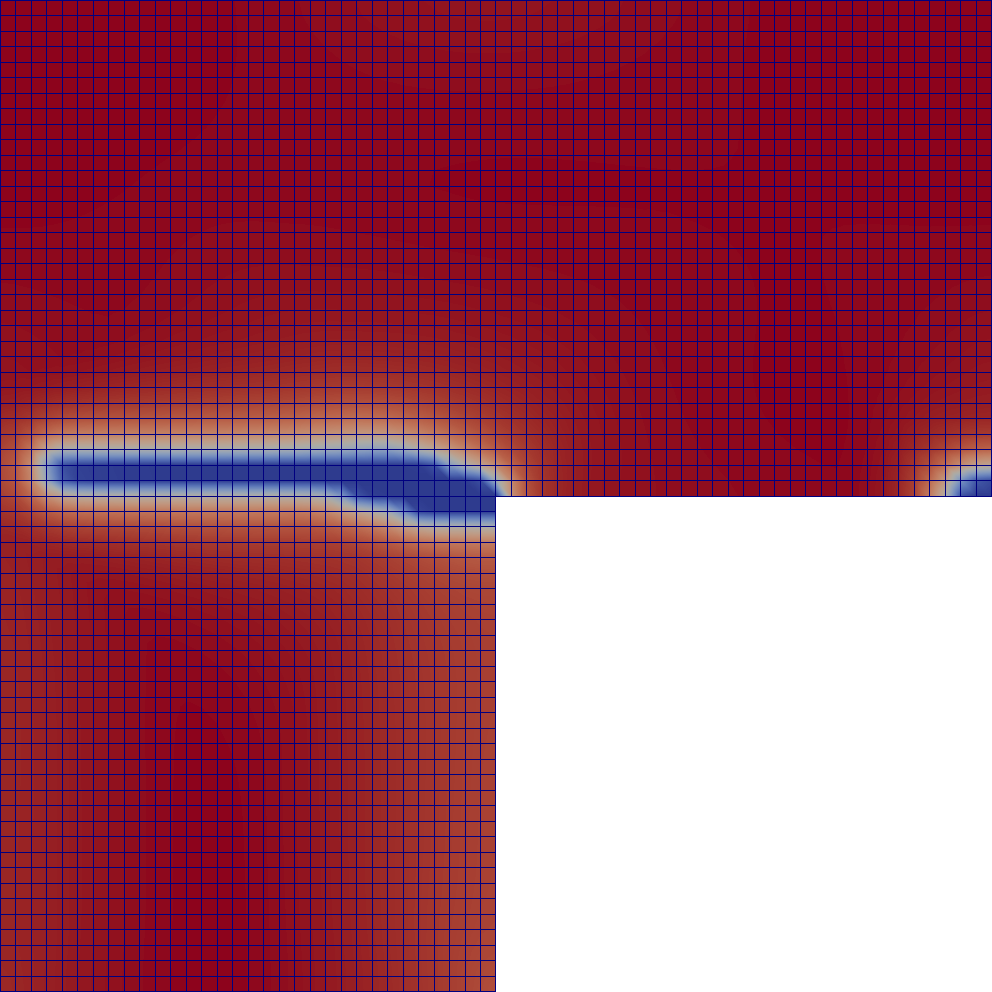}
		
		\caption{L-Shaped Panel.
			Resulting fracture pattern at times $t = 0.22s, t = 0.3s$ (top) and $t = 1.45s, t = 2s$ (bottom)}
		\label{lshape:pic}
	\end{figure}

	The computational domain $\Omega_L$ is given by an L-shape $(0, 500)^2 \setminus \left( (250, 500) \times (0, 250) \right)$, as depicted in \cref{lshape:configuration} (measures given in $mm$).
	Fracture propagation is driven by cyclic displacement boundary conditions $u_y(t)$ on $(500 - 30, 500) \times \{ 250 \}$, i.e. the lower boundary of the upper right part of the L-shape.
	The specimen is fixed on the lower part, i.e. $u = 0 \text{ on } (0, 250) \times \{0\}$.
	No additional pressure is enforced inside the crack, hence $p = 0$.
	The time-dependent displacement (given in $mm$) is defined by 
	$$u_y(t) := \begin{cases}  t & 0 s \leq t < 0.3 s, \\ 0.6 - t & 0.3 s \leq t < 0.8 s, \\ -1 + t & 0.8 s \leq t < 2 s, \end{cases}$$
	leading to a final displacement of $1 mm$ after $2 s$.
	A constant time-step of $dt := 10^{-3}$ was employed.
	Again, we chose $\varepsilon := h$ for the fracture parameter.
	The load is evaluated on the top boundary $\Gamma_{up} := (0,500) \times \{ 500 \}$.
	The full list of parameters is summarized in \cref{lshape:summary,lshape:parameters}.
	\Cref{lshape:pic} shows the resulting crack pattern at different points in time.
	
	\begin{table}[ht]
		\renewcommand{\arraystretch}{1.3}
		\centering
		\summary{data/summary_lshape_refinement_heps/summary.csv}

		\caption{L-Shaped Panel. 
			Important values and parameters for different refinement levels $l$.}
		
		\label{lshape:summary}
	\end{table} 
	
	\begin{table}[ht]
		\renewcommand{\arraystretch}{1.3}
		\centering
		\begin{tabular}{c|c|l}
			Variable & Value & Unit \\ \hline
			$\mu$ & $ 10.95 $ & $[kN/mm^2]$ \\
			$\lambda$ & $ 6.16 $ & $[kN/mm^2]$ \\
			$G_c$ & $8.9 \cdot 10^{-5}$ & $[kN / mm]$ \\
			$\kappa$ & $10^{-10}$ & $[1]$ \\
			$\varepsilon$ & $11$ & $[mm]$ \\
			$dt  \ [s]$ & $10^{-3}$ & $[s]$ \\
		\end{tabular}
		
		\caption{L-Shaped Panel.
			Parameters.}
		
		\label{lshape:parameters}
	\end{table}
	
	\subsubsection{Numerical Results}
	
	For this second test, the number of Active-Set iterations is shown for $\varepsilon = h$ and $\varepsilon = 22.0$. Random peaks happen in either case, but seem to cluster at different times. The iterations counts increase with higher levels towards the end of the simulation in case of constant $\varepsilon$. The $h$-dependent scenario is more robust in this regard.
	
	The number of GMRES iterations presented in \cref{lshape:iterations} shows similar characteristics as before.
	Again, the iteration counts increase once the fracture starts growing at $t \sim 0.3s$ and again at $t \sim 0.8s$.
	In the first half of the simulation, our linear solver is robust with respect to $h$-refinement.
	During the second half, the number of iterations increases mildly on refinement in the case of $\varepsilon = h$, but stays almost constant for $\varepsilon = const$.
	On the finest grid ($l=8$), a huge increase can be observed in either case at $t \sim 0.8s$.
	in case of constant $\varepsilon$, such irregular behavior can also be observed in the first half of the simulation.
	
	These results were obtained using the Miehe splitting, but similar findings were gathered for the Amor splitting.
	The solver also works if no splitting is applied ($\sim 2-7$ iterations), however, no crack initiation occurred in this case.
	
	\begin{figure}[ht]
		\centering
		\includegraphics[width=0.45\textwidth]{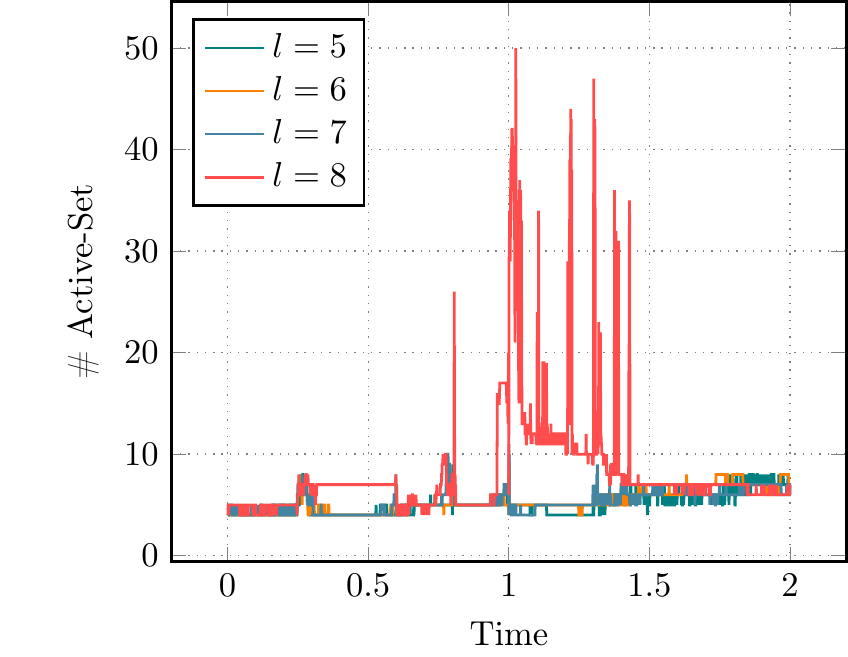}
		\includegraphics[width=0.45\textwidth]{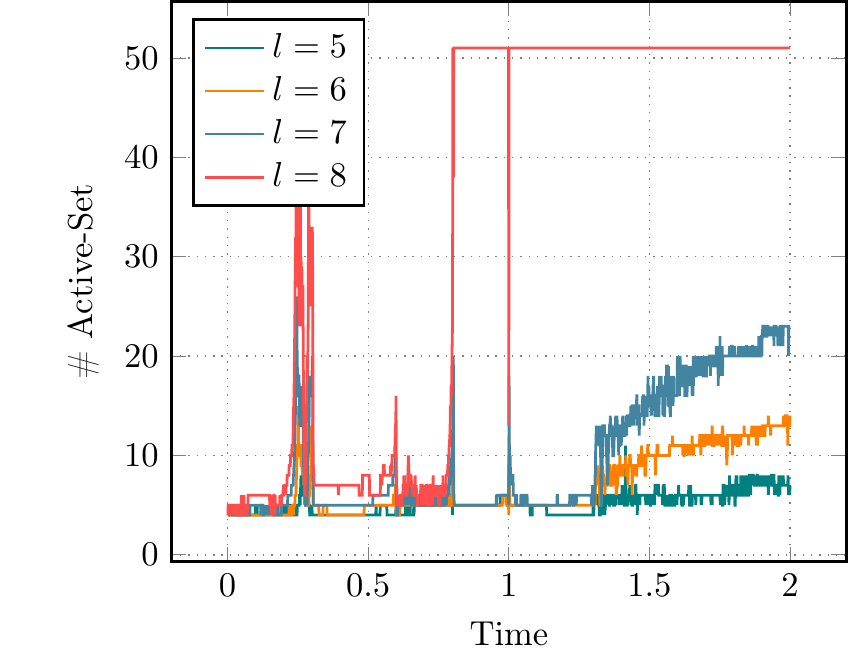}
		
		\caption{L-shaped panel.	
			Number of Active-Set steps over time for different refinement levels $l$ with $\varepsilon = h$ and $\varepsilon = 22 mm$ (right).}
		
		\label{lshape:activeset}
	\end{figure}
	
	\begin{figure}[ht]
		\centering
		\includegraphics[width=0.45\textwidth]{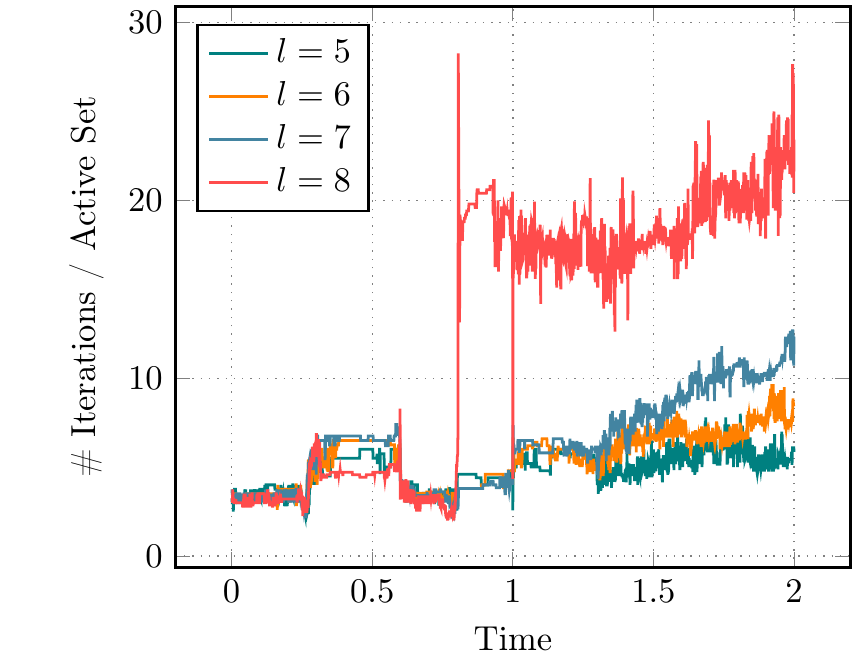}
		\includegraphics[width=0.45\textwidth]{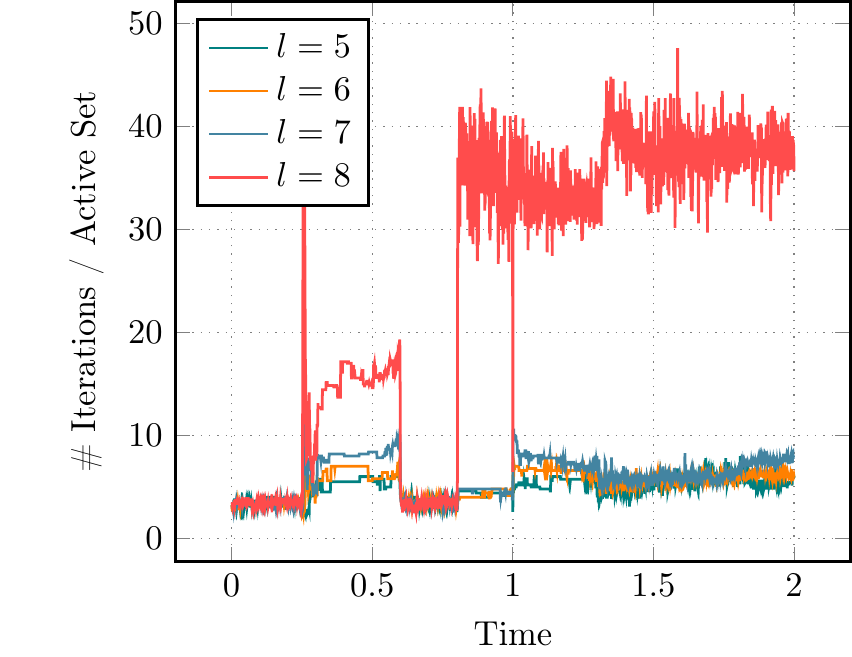}
		
		\caption{L-shaped panel.
			Iterations of the linear solver per active set step over time for different refinement levels $l$ with $\varepsilon = h$ (left) and $\varepsilon = 22 mm$ (right).}
		
		\label{lshape:iterations}
	\end{figure}
	
\begin{figure}[ht]
		\includegraphics{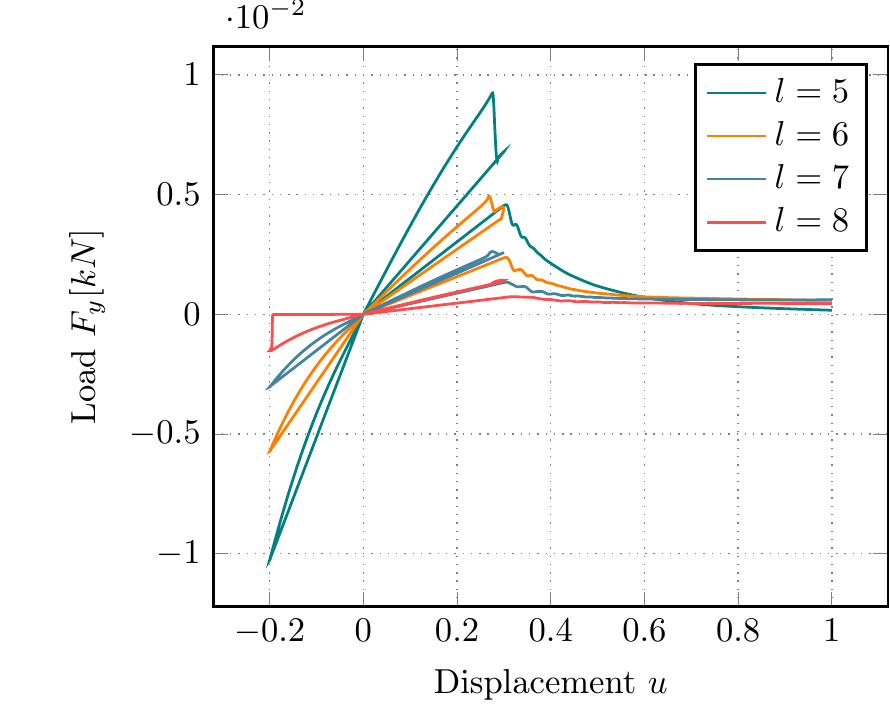}
		\includegraphics{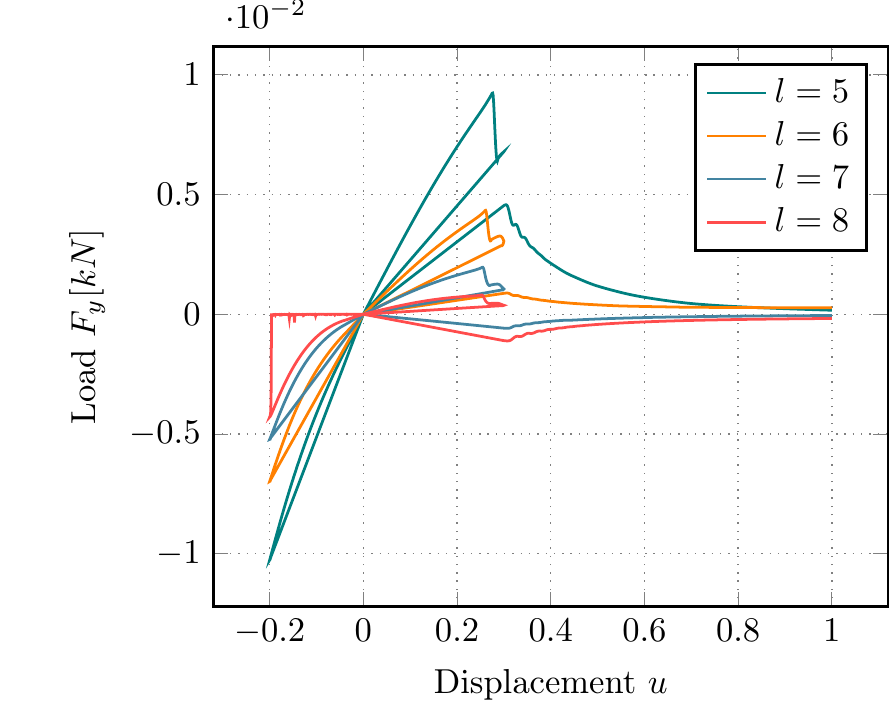}
		
		\caption{L-shaped panel.
			Load -- displacement curves for different refinement levels $l$.
			The left image shows the results for $\varepsilon = h$, whereas on the right we have a fixed value $\varepsilon = 22$.}
		
		\label{lshape:load_displacement}
	\end{figure}

	\Cref{lshape:load_displacement} shows the load-displacement curves for different refinement levels.
	Our findings look similar to results found in literature in case of $\varepsilon = h$.
	Since the phase-field parameter $\varepsilon$ changes during refinement, convergence of the loading curves is not necessarily expected.
	Nonetheless, the overall shape of the presented curves agree mostly on each other.
	A major difference happens on the finest grid ($l = 8$), where the load drops to zero at the turning point $u=-0.2m$, which does not happen on the coarser levels.
	
\begin{figure}[ht]
		\includegraphics{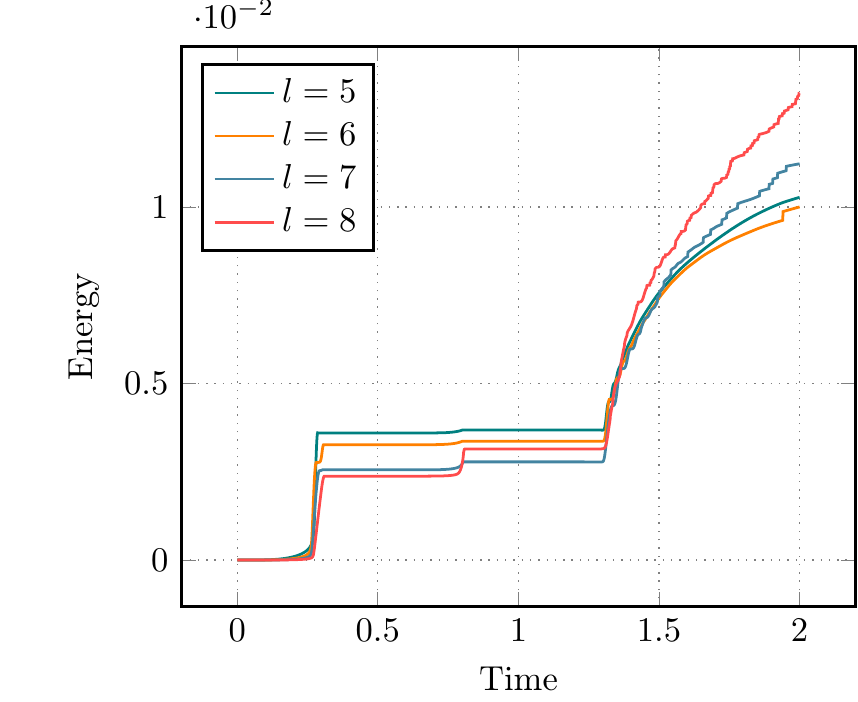}
		\includegraphics{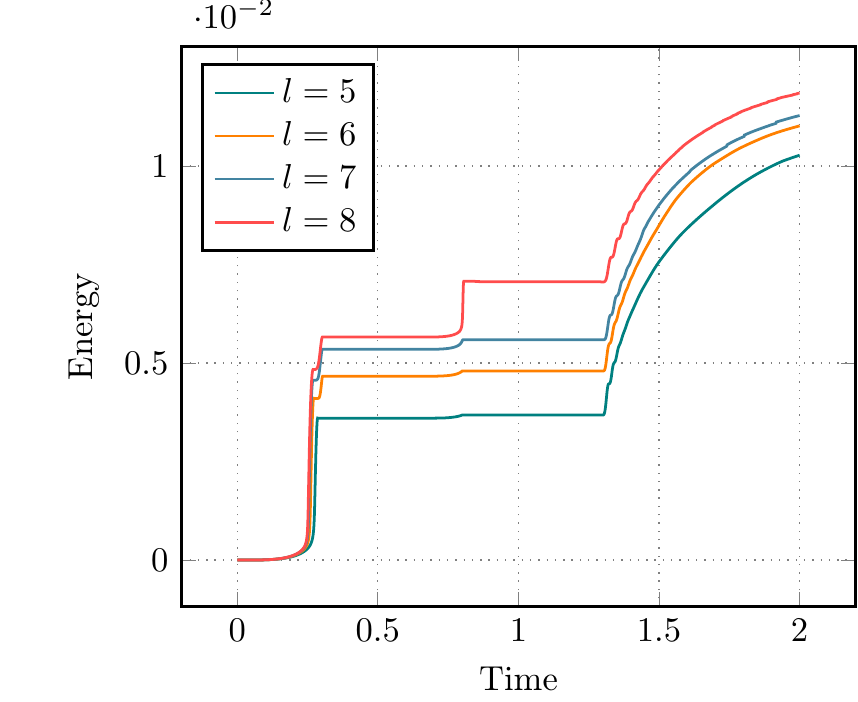}
		
		\caption{L-shaped panel.
			Crack energy for different refinement levels $l$ with $\varepsilon = h$ (left) and $\varepsilon = 22$ (right).}
		
		\label{lshape:crackenergy}
	\end{figure}

\begin{figure}[ht]
		\includegraphics{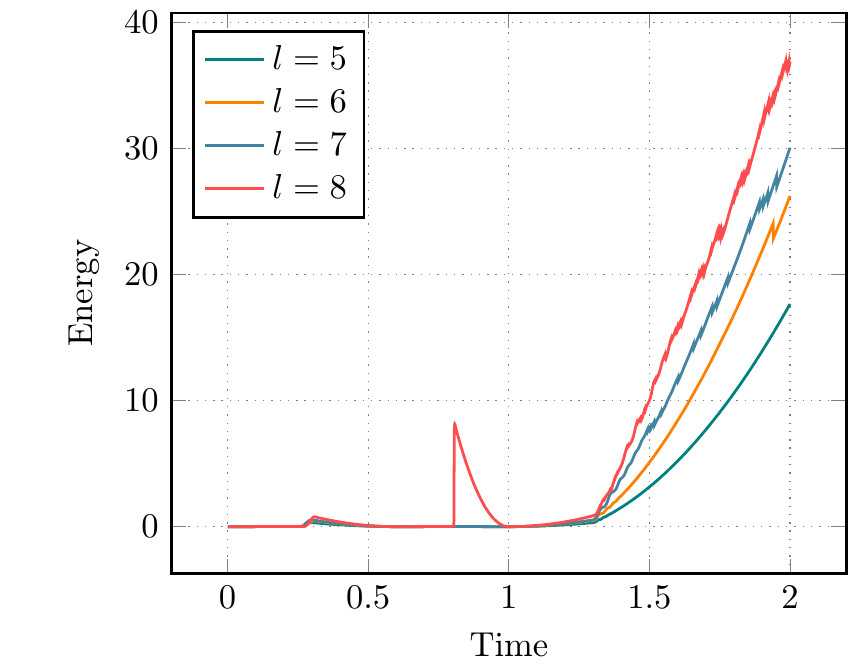}
		\includegraphics{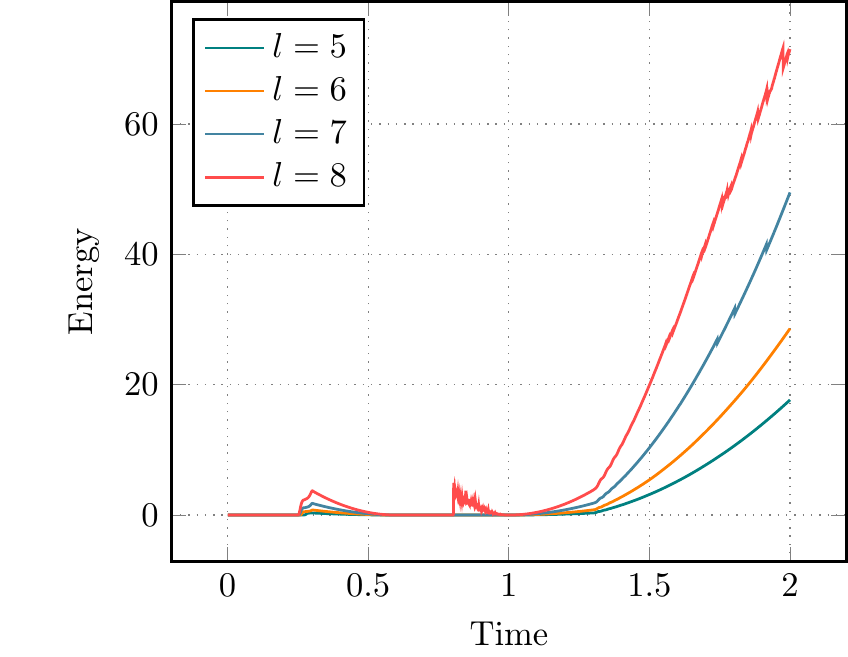}
		
		\caption{L-shaped panel.
			Bulk energy for different refinement levels $l$ with $\varepsilon = h$ (left) and $\varepsilon = 22$ (right).}
		
		\label{lshape:bulkenergy}
	\end{figure}

	The same effect is also observed in the right image of \cref{lshape:load_displacement} (constant $\varepsilon$).
	In addition, the refinements $l=7$ and $l=8$ yield negative loading values during the second cycle, which is not as expected.
	Furthermore, small irregularities are visible at turning points, i.e. step-like curves somewhere between displacements of $0.2-0.4 mm$.
	Similar artifacts are also visible in the crack and bulk-energy, see \cref{lshape:crackenergy,lshape:bulkenergy}.
	This test case seems to be very sensitive to the actual choice of
	parameters, elements and meshes, 
	which is also reported in \cite{Wi17,MaWiWo19}.
	
	\clearpage
	\subsection{$3d$ L-Shaped Panel}
	
	\subsubsection{Problem Description}
	
\begin{figure}[ht]
		\centering
		\begin{tikzpicture}
		[
		acteur/.style={
			circle,
			draw=black,
			fill=black,
			inner sep=1pt,
			minimum size=0cm
		}
		]
		\node[acteur, label={[above, xshift=4mm]\tiny $(0, 0, 0)$}] (A) at (0,0) {};
		\node[acteur, label={[above, xshift=-5mm]\tiny $(250, 0, 0)$}] (B) at (2.5, 0) {};
		\node[acteur, label=above:{\tiny $(250, 250, 0)$}] (C) at (2.5, 2.5) {};
		\node[acteur] (D) at (5, 2.5) {};
		\node[acteur, label=above:{\tiny $(500, 500, 0)$}] (E) at (5, 5) {};
		\node[acteur, label=above:{\tiny $(0, 500, 0)$}] (F) at (0, 5) {};
		\node[acteur, label={[above]\tiny $(470, 250, 0)$}] (U) at (4.4, 2.5) {};
		
		\draw (A) -- (B) -- (C) -- (D) -- (E) -- (F) -- (A);
		\draw (U.center) -- (D.center);
		\draw (A.center) -- (B.center);
		
		
		\coordinate (ddd) at (1.0, 1.20);
		
		\coordinate (UUU) at ($ (U) + (ddd) $);
		
		\node[acteur, label={[above]\tiny $(0, 500, 250)$}] (FFF) at ($ (F) + (ddd) $) {};
		\node[acteur, label={[above]\tiny $(0, 500, 250)$}] (EEE) at ($ (E) + (ddd) $) {};
		\node[acteur, label={[above]\tiny $(0, 500, 250)$}] (DDD) at ($ (D) + (ddd) $) {};
		
		\draw (D) -- (DDD);
		\draw (E) -- (EEE);
		\draw (F) -- (FFF);
		\draw[dotted] ($ (A) + (ddd) $) -- (FFF);
		\draw[dotted] (C) -- ($ (C) + (ddd) $) -- (DDD);
		\draw (B) -- ($ (B) + (ddd) $) -- ($(C)!($(B)+(ddd)$)!(D)$);
		\draw[dotted] ($ (B) + (ddd) $) -- ($ (C) + (ddd) $);
		\draw[dotted] ($ (A) + (ddd) $) -- ($ (B) + (ddd) $);
		\draw (FFF) -- (EEE) -- (DDD);

		\tikzmath{\n = 10;}
		\foreach \i in {0,...,\n}
		\draw[dotted] ($ (A) ! \i / \n ! (B) $) -- ($ ($ (A) ! \i / \n ! (B) $) + (ddd) $);
		
		\tikzmath{\m = 4;}
		\foreach \i in {0,...,\m}
		\draw[dotted] ($ (U) ! \i / \m ! (D) $) -- ($ ($ (U) ! \i / \m ! (D) $) + (ddd) $);
		
		\coordinate (M1) at ($ (D) ! 3/5 ! ($ (D) + (ddd) $) $);
		\coordinate (M2) at ($ (D) ! 4/5 ! ($ (D) + (ddd) $) $);
		\draw[<->] ($ (M1) - (0.3, 0.1) $) -- ($ (M1) - (0.3, 0) + (0, -0.75) $);
		\draw[<->] ($ (M2) - (0.3, 0.1) $) -- ($ (M2) - (0.3, 0) + (0, -0.75) $);
		\node[fill=none, label=below:{\tiny $u_y(t)$}] () at ($ (M1) ! 1/2 ! (M2) - (0, 0.35) $) {};
		
		\node[fill=none, label=center:{\tiny $u = 0$}] () at ($ (A) ! 1/2 ! ($ (B) + (ddd) $) $) {};
		
		\coordinate (O) at (4.5, 0);
		\node[fill=none, label=center:{\tiny $x$}] (X) at ($ (O) + (1, 0) $) {};
		\node[fill=none, label=center:{\tiny $y$}] (Y) at ($ (O) + (0, 1) $) {};
		\node[fill=none, label=center:{\tiny $z$}] (Z) at ($ (O) + 0.5*(ddd) $) {};
		\draw[->] (O) -- (X);
		\draw[->] (O) -- (Y);
		\draw[->] (O) -- (Z);
		
		\end{tikzpicture}
		
		\caption{$3d$ L-shaped panel. 
			Geometry and loading for the L-shaped panel test.}
		
		\label{lshape3d:configuration}
	\end{figure}
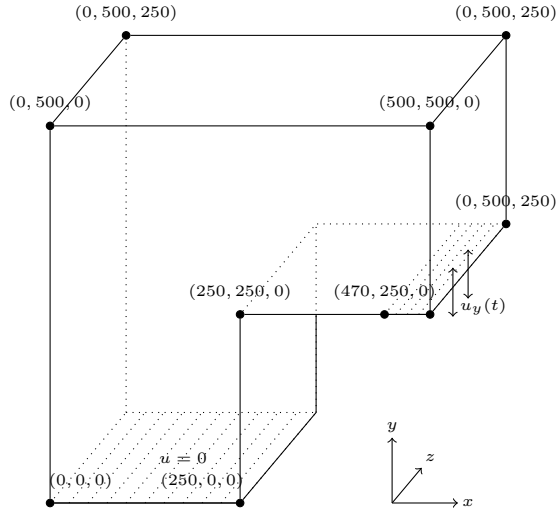

\begin{figure}[ht]
		\centering
		\includegraphics[width=0.2\textwidth]{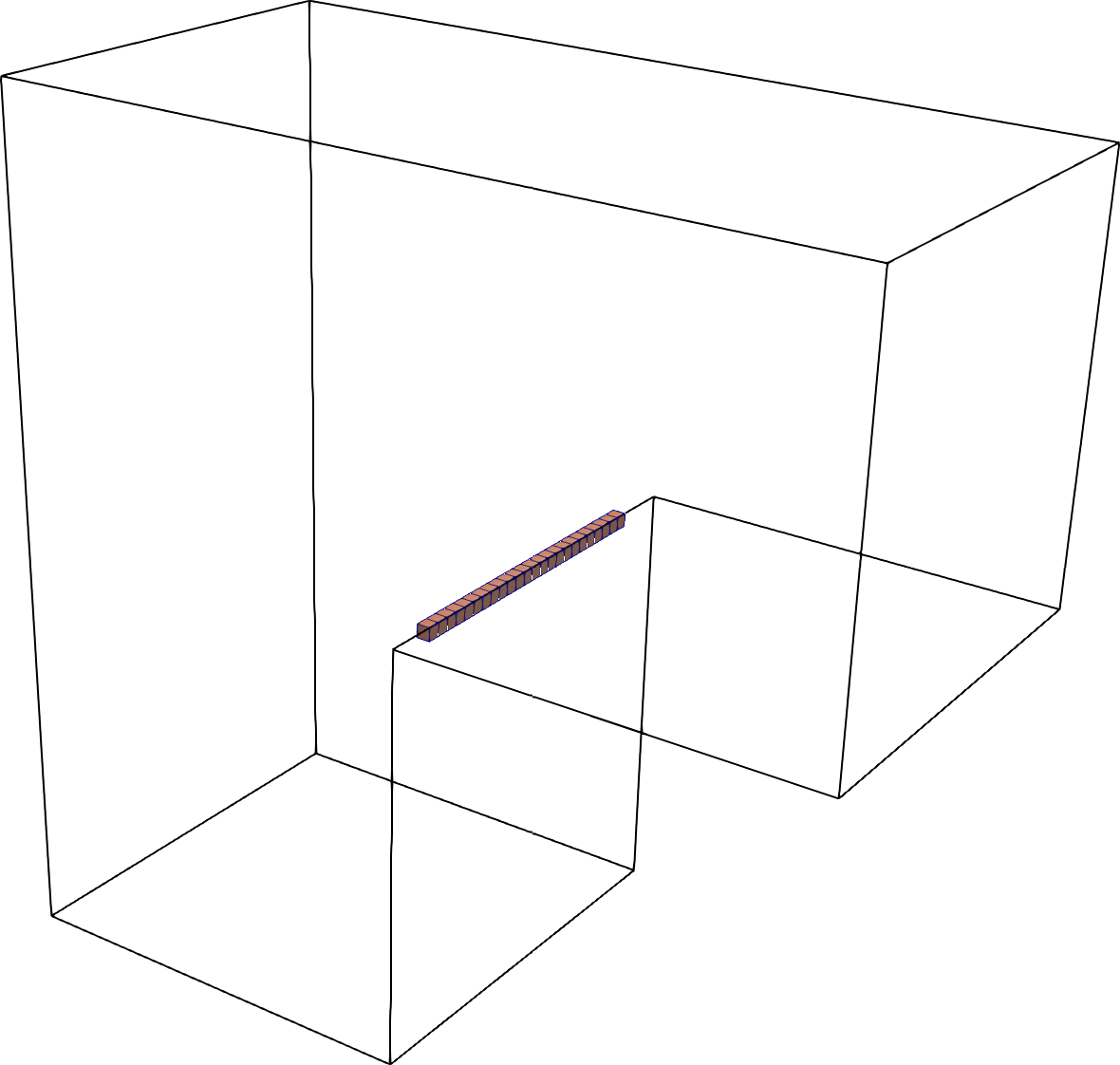}
		\includegraphics[width=0.2\textwidth]{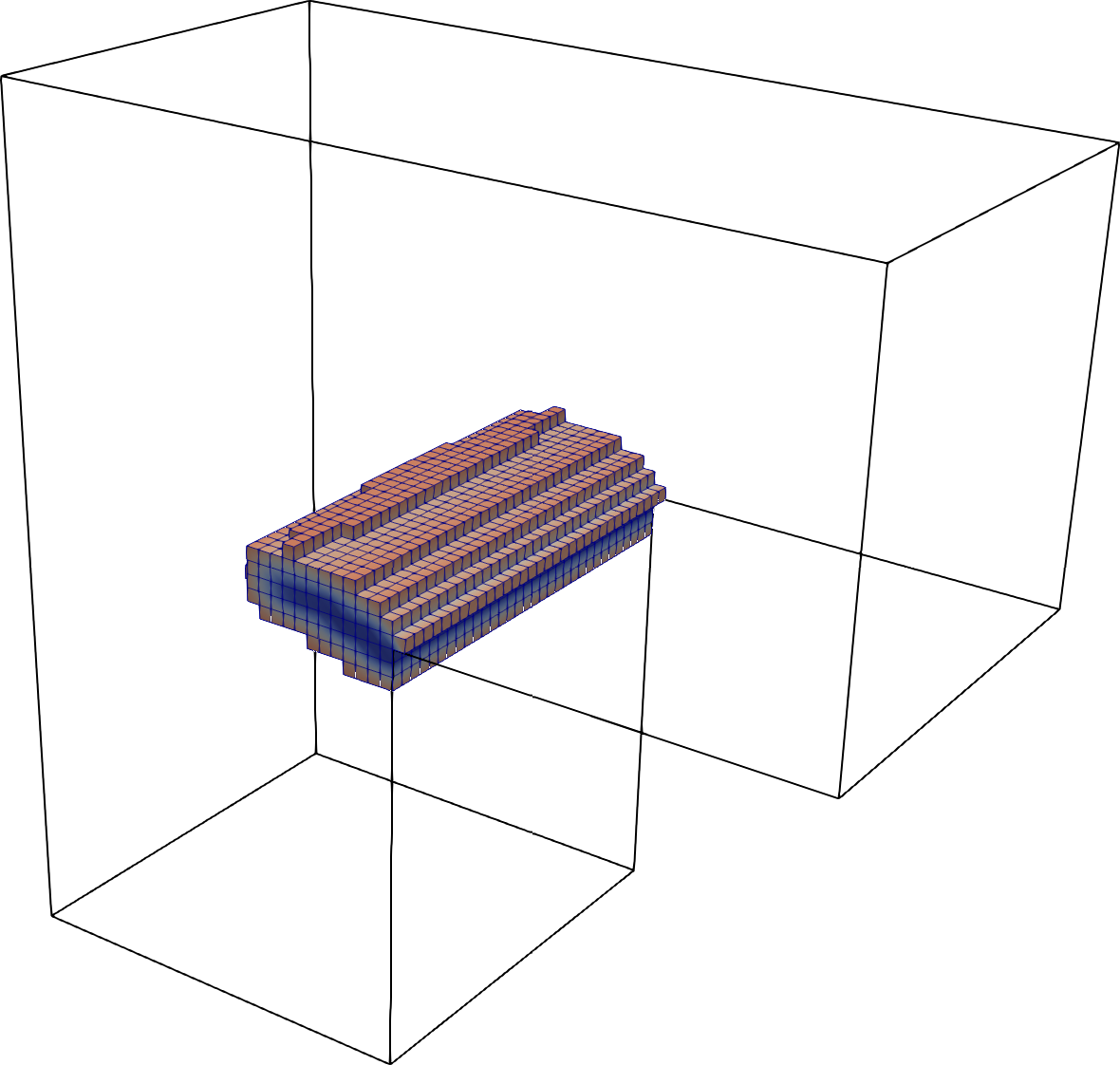}	
		\includegraphics[width=0.2\textwidth]{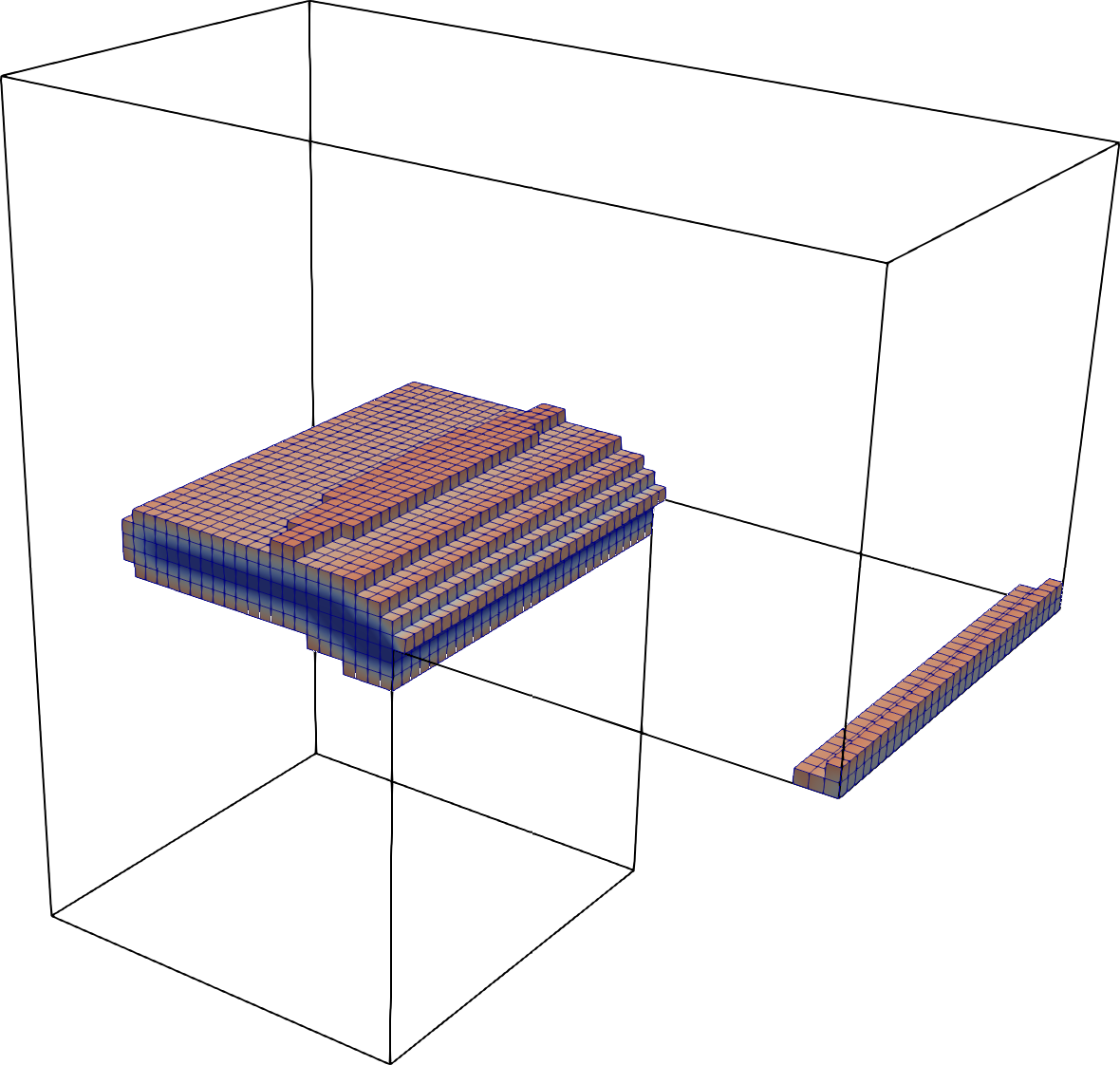}
		\includegraphics[width=0.2\textwidth]{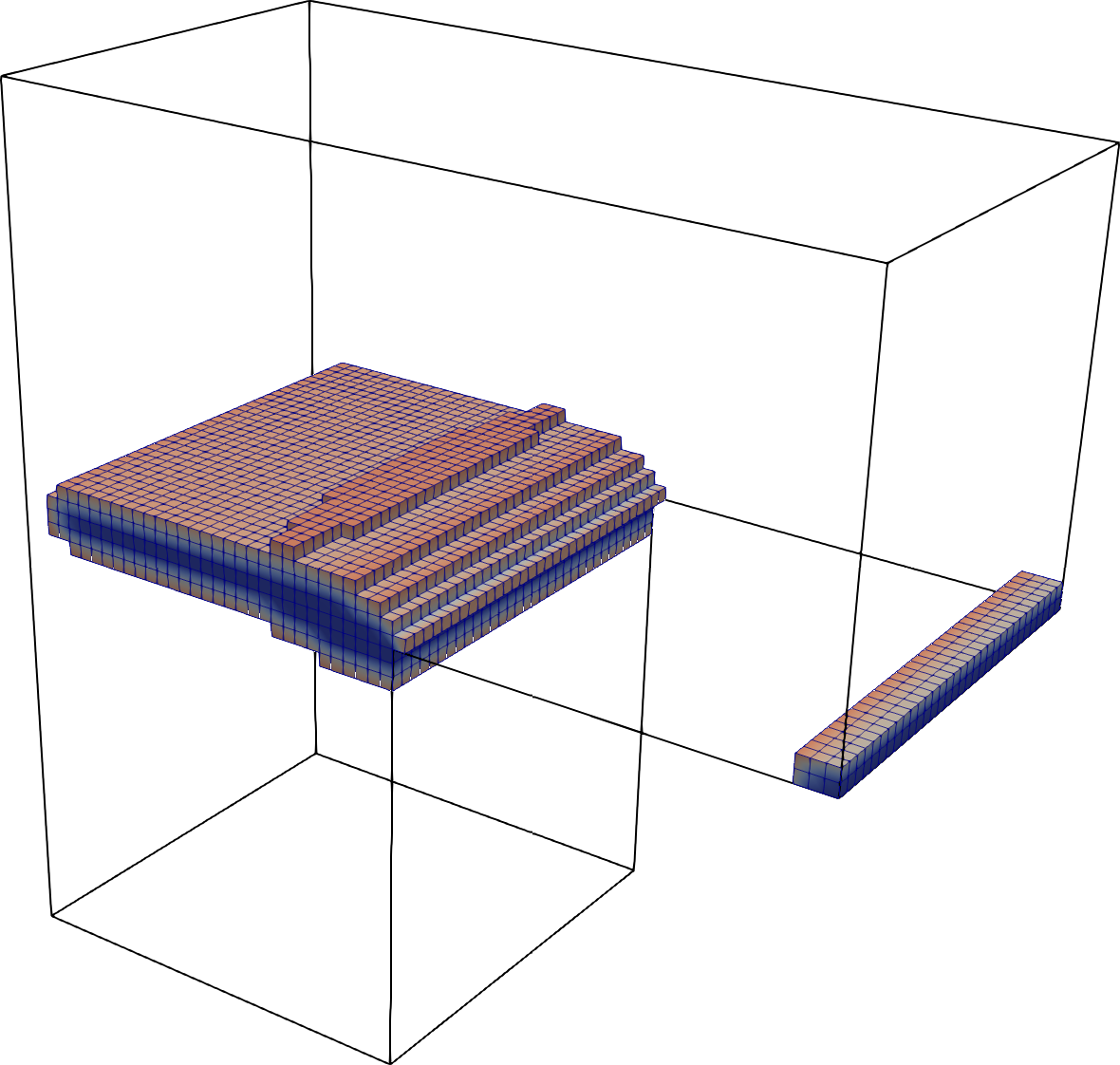}
		
		\caption{$3d$ L-Shaped Panel.
			Resulting fracture pattern at times $t = 0.22s, t = 0.3s$ (top) and $t = 1.45s, t = 2s$ (bottom). Only parts with $\varphi < 0.75$ are shown.}
		\label{lshape3d:pic}
	\end{figure}

	In this section, we consider a $3d$ version of the previous L-shaped
	panel test. A related numerical study is \cite{MeBoKh15}. A related
	experimental test was previously proposed in \cite{Wi01}. 
	The computational domain is given by extruding $\Omega_L$ from before into z-direction.
	Hence, $\Omega_L^{3d} := \Omega_L \times (0, 250)$, as depicted in \cref{lshape3d:configuration} (measures given in $mm$).
	The same cyclic displacement boundary conditions are applied as before on the extruded area $(500 - 30, 500) \times \{ 250 \} \times (0, 250)$.
	The specimen is fixed on the bottom, i.e. $u = 0 \text{ on } (0, 250) \times \{0\} \times (0, 250)$.
	Again, the loading force is evaluated on the upper boundary $\Gamma_{up} := (0,500) \times \{ 500 \} \times (0, 250)$.
	The full list of parameters is summarized in \cref{lshape3d:summary,lshape3d:parameters}.
	\Cref{lshape3d:pic} illustrates the resulting fracture.
	
	\begin{table}[ht]
		\renewcommand{\arraystretch}{1.3}
		\centering
		\summary{data/summary_lshape3d/summary.csv}
		
		\caption{$3d$ L-Shaped Panel. 
			Important values and parameters for different refinement levels $l$.}
		
		\label{lshape3d:summary}
	\end{table} 
	
	\begin{table}[ht]
		\renewcommand{\arraystretch}{1.3}
		\centering
		\begin{tabular}{c|c|l}
			Variable & Value & Unit \\ \hline
			$\mu$ & $ 10.95 $ & $[kN/mm^2]$ \\
			$\lambda$ & $ 6.16 $ & $[kN/mm^2]$ \\
			$G_c$ & $8.9 \cdot 10^{-5}$ & $[kN / mm]$ \\
			$\kappa$ & $10^{-10}$ & $[1]$ \\
			$\varepsilon$ & $11$ & $[mm]$ \\
			$dt  \ [s]$ & $10^{-3}$ & $[s]$ \\
		\end{tabular}
		
		\caption{$3d$ L-Shaped Panel.
			Parameters.}
		
		\label{lshape3d:parameters}
	\end{table}
	
	\subsubsection{Numerical Results}
	
\begin{figure}[ht]
		\centering
		\includegraphics[width=0.45\textwidth]{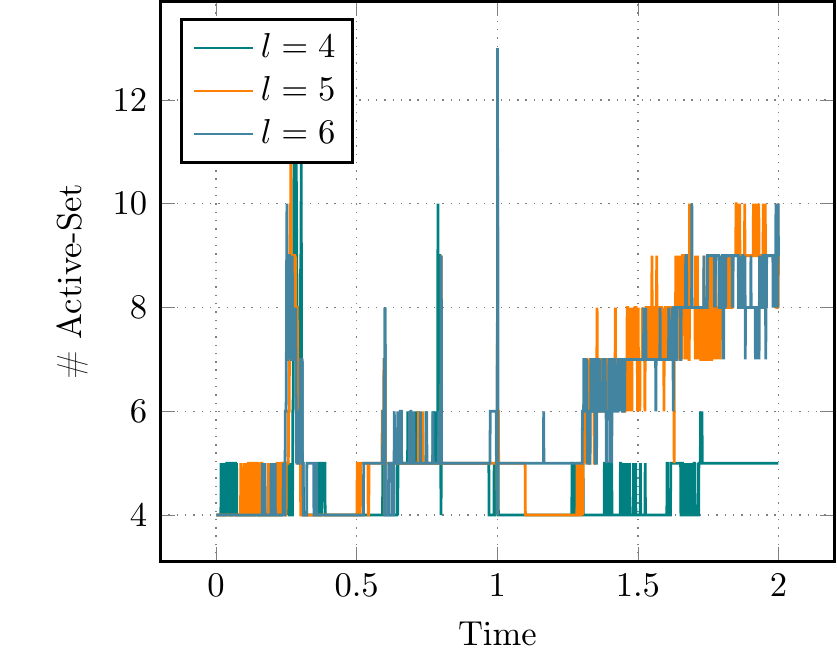}
		\includegraphics[width=0.45\textwidth]{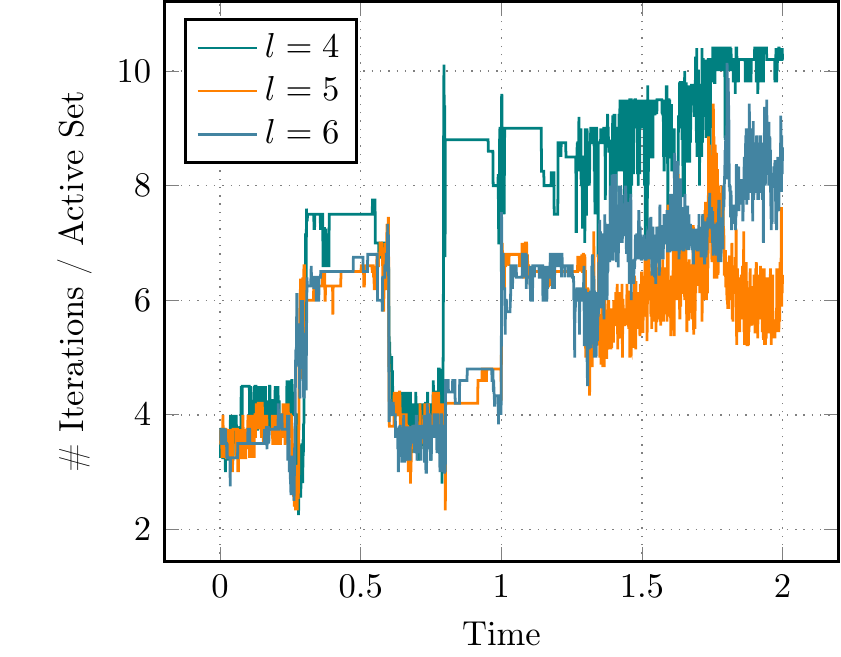}
		
		\caption{$3d$ L-shaped panel.
			Left: number of Active-Set iterations over time.
			Right: number of iterations of the linear solver per active set step over time for different refinement levels $l$ with $\varepsilon = h$.}
		
		\label{lshape3d:iterations_active}
	\end{figure}

	The results for the $3d$ extension of the L-shaped panel test are very similar to the previously shown $2d$ setup.
	The number of Active-Set iterations is in the same range as before, but less peaks are visible in the $3d$ setting, see \Cref{lshape3d:iterations_active} (left).
	\Cref{lshape3d:iterations_active} (right) shows that the number of GMRES iterations per Active-Set step is in again very similar, i.e. roughly $4-10$ iterations.
	
	We also observe a similarly increased iteration count at times when the fracture is growing.
	
\begin{figure}[ht]
		\includegraphics{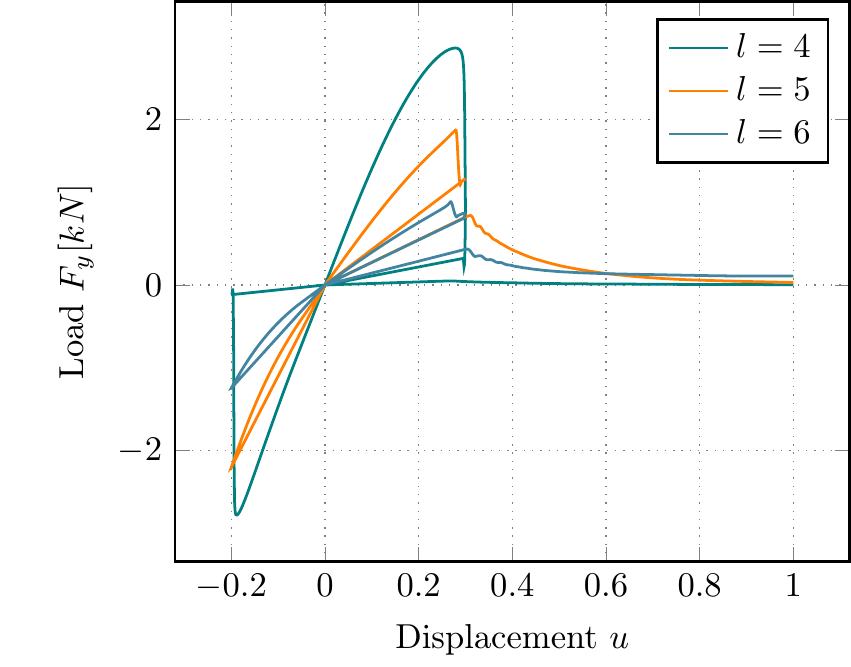}
		\includegraphics{plot/auto_lshape_refinement_heps_load_displacement}
		
		\caption{$3d$ L-shaped panel.
			Load -- displacement curves for different refinement levels $l$ in the $3d$ case (left).
			For comparison, we show the respective results for $2d$ case again (right).}
		
		\label{lshape3d:load_displacement}
	\end{figure}

\begin{figure}[ht]
		\includegraphics{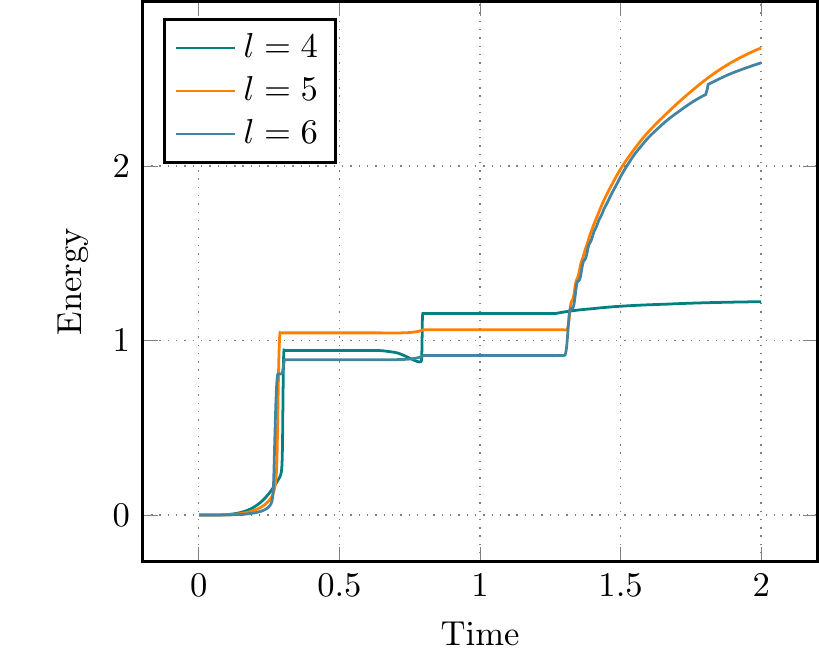}
		\includegraphics{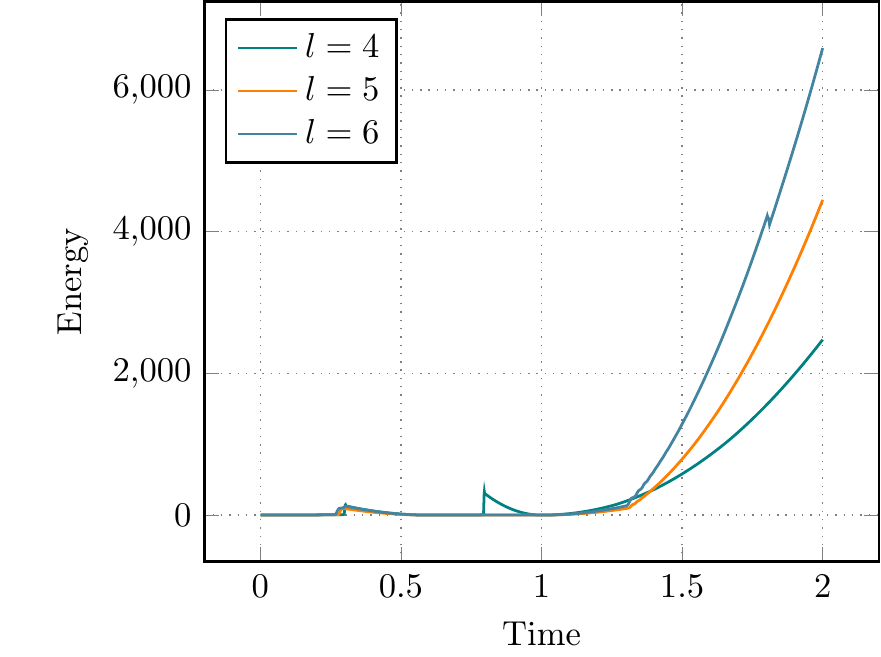}
		
		\caption{$3d$ L-shaped panel.
			Crack and bulk energy for different refinement levels $l$ with $\varepsilon = h$.}
		
		\label{lshape3d:energy}
	\end{figure}

	Also, the shape of energy and loading curves (\cref{lshape3d:load_displacement,lshape3d:energy}) look very similar in the $2d$ and $3d$ case.
	The exception being the coarsest $3d$ grid, which shows a very roundly shaped loading curve, which does not agree with the $2d$ setting.
	However, the numbers do not even remotely agree, as seen in \cref{lshape3d:load_displacement,lshape3d:energy}.
	
	\clearpage
	\section{Conclusions and Future Work}
	
	In this work, we developed a matrix-free geometric multigrid solver
	for fracture propagation problems
	using a phase-field description.
	Our numerical experiments show that this solver can handle a variety of test cases within a reasonable number of iterations.
	Although we have only shown results for the AT 2 model, it shall be noted that the solver also works for the simpler AT 1 model, with iteration counts usually below the presented values.
	Limitations of this method appear primary in case of stress-splittings.
	The improvement of the solver in these cases is part of our future work.
	In order to deal with larger simulations (in particular in 3d), incorporating adaptivity and/or parallelization is necessary; both of which are part of our ongoing work.
	Another point of concern is the Active-Set method.
	In some simulations, convergence of this non-linear solver is observed
	to be very slow as also indicated in fundamental studies on
	primal-dual active set methods, e.g., \cite{Ha15,CuHaRo15}.
	Hence, it may be beneficial to further investigate properties of the
	active set method with the goal of further improvements.
	Ultimately, one could also think of combining the non-linear and linear solver in terms of nested iterations or non-linear multigrid methods (see e.g. \cite{Kr06,Ha94}).
	
	
	\section{Acknowledgments}
	
	This work has been supported by the Austrian Science Fund (FWF) grant P29181 `Goal-Oriented Error Control for Phase-Field Fracture Coupled to Multiphysics Problems'.
	
	
	\printbibliography
	
\end{document}